\documentclass[10pt,leqno]{amsart}
\usepackage{amsmath,amsthm,amssymb,enumitem,color,fancyvrb,prooftree,forest,abraces,tikz}
\usepackage{bm}
\usepackage[colorlinks=false,urlcolor=blue,pagebackref=true]{hyperref}
\usepackage{multicol}
\usepackage{setspace}
\allowdisplaybreaks
\theoremstyle{plain}
\newtheorem{theorem}{Theorem}[section]
\newtheorem{proposition}[theorem]{Proposition}
\newtheorem{lemma}[theorem]{Lemma}
\newtheorem{corollary}[theorem]{Corollary}

\newtheorem{claim}{Claim}
\theoremstyle{definition}
\newtheorem{definition}[theorem]{Definition}
\newtheorem{notation}[theorem]{Notation}
\newtheorem{example}[theorem]{Example}

\newtheorem{remark}[theorem]{Remark}
\newtheorem{note}[theorem]{Note}

\newtheorem{problem}[theorem]{Problem}
\newtheorem{convention}[theorem]{Convention}

\newcommand{\innerthmname}{}% initialize

\theoremstyle{definition}

\newcommand{\N}{\mathbb{N}}
\newcommand{\vp}{\varphi}
\newcommand{\E}{\bigvee}
\newcommand{\A}{\bigwedge}

\newcommand\wse{winning strategy for Eloisa }
\newcommand{\guess}{\mathtt{guess}\,{:}\  }
\newcommand{\query}{\mathtt{query}\,{:}\ }
\newcommand{\reply}{\mathtt{reply}\,{:}\ }
\newcommand\pGamma{\dot{p}}
\newcommand{\cut}{\mathrm{cut}}
\newcommand{\DR}{\mathsf{DR}}
\newcommand{\on}{\mathbb{O}}
\newcommand{\Ord}{\mathrm{ON}}
\DeclareMathOperator{\dom}{\mathit{dom}}

\DeclareMathOperator{\height}{\mathit{ht}}
\newcommand{\dt}{\mathit{dh}}     %depth

\newcommand{\conc}{{{}^\smallfrown}}
\newcommand{\restr}{\upharpoonright}
\newcommand{\lr}{\langle\rangle}

\newcommand{\subsetm}{\subset_{m}}
\newcommand{\hh}{h}   % height function
\newcommand{\party}{\imath}

\def\ru#1#2#3{\prooftree #1 \justifies #2 \using{\text{#3}} \endprooftree}

\def\pa{{\sf PA}}
\DeclareMathOperator{\TI}{\mathsf{TI}}
\newcommand{\imp}{\rightarrow}

\title[Peano Arithmetic, games and descent recursion]{Peano Arithmetic, games and descent recursion}
\thanks{The author's research was supported by the Alexander von Humboldt foundation.}

\author[E.\ Frittaion]{Emanuele Frittaion}
\address{Department of Mathematics, Technische Universit\"{a}t Darmstadt, Germany}
\begin{document}
	
	
\begin{abstract} 
We analyze Coquand's game-theoretic interpretation of Peano Arithmetic \cite{C95} through the lens of  elementary descent recursion \cite{FS95}. In Coquand's game semantics, winning strategies correspond to infinitary cut-free proofs and  cut elimination corresponds to   {\em debates} between these winning strategies.  The proof of cut elimination, i.e., the proof that such debates eventually terminate,  is by  transfinite induction on certain {\em interaction} sequences of ordinals.    In this paper, we provide a direct implementation of Coquand's proof, one that allows us to describe winning strategies by descent recursive functions. As a byproduct, we obtain yet another proof of well-known results about provably recursive functions and functionals. 
\end{abstract}
	
\maketitle

\section{Coquand's game semantics}

In \cite{C95} (see also \cite{BBC98,C91,C97}) Coquand presents a constructive game semantics  for classical arithmetic, building on Novikoff's notion of regularity \cite{N43}. Regularity  corresponds to cut-free provability   in a suitable infinitary sequent calculus. The basic idea is that an infinite (well-founded) cut-free  proof  of a sentence $\vp$ can be seen as a winning strategy for one player in a certain two-player game $G(\vp)$, with players traditionally named Eloisa for $\exists$ and Abelard for $\forall$.  The paths of the proof tree represent all possible plays that agree with Eloisa's strategy. Since every true sentence of arithmetic has an  elementary recursive cut-free proof, there are always elementary recursive winning strategies for Eloisa. In this sense, Coquand's semantics is constructive.

By contrast, Tarski game semantics is not.  Recall that in a  Tarski game for a sentence $\vp$,  Eloisa and Abelard must play consecutive immediate subformulas  of $\vp$ until  a literal is obtained, i.e., an atomic formula or a negation thereof  (we can assume that the sentence is written by using $\neg, \lor, \land, \exists, \forall$ and $\neg$ only occurs in front of atomic formulas).  
On each turn, Eloisa plays  a subformula $\psi_i$, resp.\ $\psi(\bar n)$, whenever the last move is $\psi_0\lor\psi_1$, resp.\ $\exists x\, \psi(x)$.  Dually, Abelard plays a subformula $\psi_i$, resp.\  $\psi(\bar n)$, whenever the last move is $\psi_0\land\psi_1$, resp.\  $\forall x\, \psi(x)$. Eloisa wins if the last move on the table is a true literal.  Any such game  ends in  finitely many steps (the rank of $\vp$). Even though  either player has a winning strategy --   Eloisa  if $\vp$ is true, and Abelard otherwise  --   Eloisa need  not always have a recursive one.  In fact, there are sentences that are valid in classical logic  and have no recursive winning strategy for Eloisa.  A standard example is
\[    \forall x\, \exists y\, (\exists z\, R(x,z)\imp  R(x,y)), \]
where $R(x,y)$ is recursive but  $\exists y\, R(x,y)$ is not. Ironically, Tarski game semantics fails to be constructive because winning strategies correspond to intuitionistic cut-free proofs. 
On the other hand, in a Coquand   game  Eloisa can always backtrack and make a new guess as often as she pleases. Essentially a greedy search will do. But we can do better if the sentence $\vp$ is a theorem of $\pa$ (Peano Arithmetic). 

The real content of \cite{C95} is a game-theoretic proof of cut elimination. In this context, cut elimination takes the form of a {\em debate} between winning strategies (i.e., the cut-free proofs of the left and right premise of the cut). Debates are characterized  by certain sequences of integers, which he calls {\em interaction sequences}, and the proof that every debate terminates is based on the fact  that any sequence of ordinals $\alpha_0\alpha_1\alpha_2\cdots$ governed by an interaction sequence must necessarily be finite. In this paper, we aim to proof-mine Coquand's termination argument and draw some proof-theoretic consequences for $\pa$. In keeping with the analogy between strategies and proofs,   we want to show how to compute the steps of a proof along with ordinals witnessing that the proof tree is well-founded. To do so, we will deploy the framework of elementary descent recursion \cite{FS95}.  Descent recursive functions  provide an elegant characterization of  $\alpha$-recursion in the sense of proof theory (see, e.g., \cite[Chapter 4]{SW12}) and can be intuitively understood as a   bounded search  to be executed over an  ordinal number of steps, rather than just a finite number of steps.\footnote{For an in-depth  analysis of descent recursion below $\varepsilon_0$, including a tight comparison with the Hardy hierarchy, we refer the reader to \cite{FS95}.} It turns out that every debate associated to Coquand's cut elimination can be easily described by an elementary descent recursion.  Once we have a suitable notion of {\em descent recursive} strategy, all that remains is to prove that Coquand's cut elimination   preserves such strategies. 

Our  approach to Coquand's game-theoretic interpretation  is parallel to that of \cite{B91,FS95, M78,S77}. The key difference is that these works are all based on Sch\"{u}tte's   infinitary cut elimination for $\omega$-logic \cite{KS51}. A distinct feature of their approach is the use of specific delaying mechanisms  to handle the complexity of the {\em standard cut elimination} operator. In particular, the so-called repetition rule introduced by Mints \cite{M78}  enables to convert a top-down implementation of cut elimination (a bar recursion) into a  bottom-up procedure (a primitive recursion). 
More in detail, Mints' primitive recursive operator  $E$   performs full cut elimination on locally correct proofs.  It is defined by iterating the operator R, which lowers the cut rank of a proof by one and it does so in an elementary recursive way.  Notice that Mints' operator is defined  everywhere and so it also applies to proofs that are not well-founded.
By contrast, Coquand's  procedure  removes a single cut -- the bottom-most one -- at once, sidestepping   the {\em standard reduction} of a cut into smaller ones (the reduction simulated by  $R$). The procedure gives rise to a  partial recursive operator that   only  terminates on well-founded proofs.  One could mirror the repetition  rule by using {\em dummy} moves and thus produce a  total primitive recursive variant of Coquand's operator. However,  we have no reason to do so. The delaying effects of cut elimination are simply hidden in the computation, resulting in longer descent recursions. Finally, we use descent recursion as  in \cite{FS95}. The difference is that they  consider descent recursive proofs (strategies) with cuts and repetitions.  In particular,  they  need to show that such proofs  are preserved under the operator $R$, which   calls for a rather delicate treatment of higher type functionals  (see \cite[Section 1.1]{FS95}). In a way, our approach is more direct and cost-effective.    

We point out that  Coquand's cut elimination produces an  increase of the height by a tower of $\nu$-many exponentials, where the number $\nu$ is the depth (changes in polarity) of the cut formula.  In the case of  debates between finite strategies,   much better  bounds haven been  obtained. In fact, Aschieri \cite{A17} showed that the height of the tower depends on the minimum of the so-called backtracking levels of the winning strategies, which are  independent  of the depth $\nu$ and in the worst case scenario are  equal to  $\nu-2$ for $\nu\geq 2$.   It would be interesting to lift such results to classical arithmetic and to compute ordinal bounds that only depend on the backtracking level of the winning strategies. Our purpose here, however, is a bit different. We are not interested in the bounds   {\em per se}, but in how to implement Coquand's cut elimination in order to obtain descent recursive  strategies.

We will finally show how to use    strategies to easily compute   the provably recursive functions   and  the Kreisel's no-counter-example interpretation functionals of $\pa$ as well as extensions of $\pa$  by transfinite induction over elementary recursive ordinal representation systems.  Pertaining to the second application, the reader may be interested in   Tait's game-theoretic explanation \cite{T05}   of Gentzen's 1936 consistency proof.\footnote{To be precise, the posthumously published  early draft of Gentzen's first consistency proof \cite[pp.\ 201--213]{G69}.} In fact, Tait expands on  G\"{o}del's 1938  
reading  of Gentzen in terms of the no-counter-example interpretation.\footnote{See  \cite[\S VI, pp. 103--111, paragraphs 15--19]{G95}.}

There is a certain degree of freedom in  setting up a cut-free proof system whose provable formulas are exactly the regular ones. Since there is nothing special about regularity, we will adapt  Coquand's game semantics to suit the specifics of Tait calculus (see \cite{T68}) and forget about regularity altogether.

\begin{note}
The correspondence between debates and  cut elimination  in the case of Novikoff infinitary calculus (but it should translate to other calculi) has been established in \cite{H95}. A more general framework to deal with  debates and interaction sequences  is provided by  Hyland-Ong games \cite{HO00}. Among the applications, one notable example is a proof of normalization for typed $\lambda$-calculus \cite{C11}. The literature on game semantics is remarkably vast. We refer the interested reader to \cite{A17} for more information. 
\end{note}

\subsection{Structure of paper} In Section \ref{Tait} we describe in broad strokes Tait (infinitary) calculus and in Section \ref{Peano} we outline a (finitary) Tait calculus for Peano Arithmetic; at a first reading,  these two sections
should perhaps be referred to only as necessary.  Section \ref{game} introduces a game semantics directly inspired by (the cut-free segment of) Tait calculus; in particular,  we fix some useful terminology to describe formula trees and introduce a  game-theoretic notion related to such tree  representations that will be crucial later on. The heart of the paper are Sections \ref{debates},  \ref{finite depth debates} and \ref{computing}: the first   introduces  the notion of debate in general, the second  discusses  debates of finite depth, and the third provides an explicit construction of a debate.  Section \ref{descent recursion} introduces the main tool to analyze the complexity of winning strategies. Finally, in Section \ref{peano games}, we put it all together and prove our main results.

\section{Tait  calculus} \label{Tait}

For the reader's convenience, we briefly recall the set-up of Tait calculus and assume familiarity with e.g.\   \cite{T68}.  One may safely take the game-theoretic view of this calculus  at face value (Section \ref{game}).  
\begin{definition}[formulas]
Let $P$ be a collection of {\em atoms}. Formulas (in negation normal form) are inductively generated by  the following rules:
\begin{itemize}[leftmargin=5mm]
\item literals of the form $p$ and $\neg p$ with $p\in P$ are formulas; 	
\item if $\vp_i$ is a formula for every $i\in I$ with $I\neq\emptyset$, then $\E_{i\in I} \vp_i$ ($\E$-formulas) and  $\A_{i\in I} \vp_i$ ($\A$-formulas) are formulas.	
\end{itemize}
Negation is defined by letting:
\[  \neg \neg p\simeq p\qquad \neg \E\vp_i \simeq \A\neg\vp_i\qquad \neg\A\vp_i\simeq \E\neg\vp_i.\]
\end{definition}

\begin{remark}
For our purposes, it is more convenient to deal with possibly infinite but non-empty disjunctions and conjunctions. 
\end{remark}

\begin{definition}
Sequents $\Gamma, \Delta,\ldots$ are finite sets of formulas. $\Gamma,\vp$ stands for $\Gamma\cup\{\vp\}$. 
\end{definition}

\begin{definition}[Tait calculus]
\[ \begin{array}{lcl}
& \Gamma \text{ axiom} & \\
\ru{\Gamma, \vp_i}{\Gamma,\E \vp_i}{$\E$-rule} &  & \ru{ \Gamma, \vp_{i} \ \  (i)}{\Gamma,\A\vp_i}{$\A$-rule} \\
& \ru{\Gamma,\neg \vp\ \ \ \Gamma,\vp }{\Gamma}{cut} &
\end{array}
\]
\end{definition}

A sequent $\Gamma$ is an {\em axiom} if $\Gamma\cap S\neq\emptyset$ for some fixed collection $S$ of (true) literals. Namely, it is required that for every  $p$,   either $p$ or $\neg p$ is in $S$  but not both.\footnote{It is possible to be  more general for cut elimination to go through. This is nicely done in \cite{T68}.} An infinitary Tait calculus for Peano Arithmetic is obtained by taking as $S$ the set of all true closed literals.

\section{Peano Arithmetic} \label{Peano}
The language of $\pa$ consists of $=$, $0$, $S$, together with a collection of  function and predicate symbols for primitive recursive functions and relations. We will be more specific later on.

\begin{definition}
Formulas are built up from literals (i.e., atomic and negated atomic formulas) by means of 
 $\lor$, $\land$, $\exists x$, $\forall x$.
\end{definition} 
\begin{definition}[Tait calculus for $\pa$]

\[ \begin{array}{rcr}
& 	\Gamma \text{ axiom} & \\
\ru{\Gamma, \vp\   (\text{or }\Gamma,\psi)}{\Gamma,\vp \lor \psi}{$\lor$-rule} &  & \ru{ \Gamma, \vp \ \ \  \Gamma, \psi  }{\Gamma,\vp\land\psi}{$\land$-rule} \\[5mm] 
\ru{\Gamma, \vp_x(e)}{\Gamma,\exists x\, \vp(x)}{$\exists$-rule} &  & \ru{\Gamma, \vp_x(y) }{\Gamma,\forall x\, \vp(x)}{$\forall$-rule}\\
& \ru{\Gamma,\neg \vp\ \ \ \Gamma,\vp }{\Gamma}{cut} & 
\end{array}
\]
\end{definition}

The usual restrictions on variables and terms apply ($e$ denotes a term). The axioms consist of:
\begin{itemize}[leftmargin=5mm]
	\item a collection of basic axioms; 
	\item induction axioms of the form \[ \Gamma, \neg\vp(0), \exists x\, (\vp(x)\land \neg\vp(Sx)), \forall x\, \vp(x). \]
\end{itemize}

By a basic axiom we mean a sequent $\Gamma$ such that every substitution of the free variables with closed terms (equivalently, numerals) makes at least one literal of $\Gamma$ true.  Note that the open axioms of Peano Arithmetic and the equality axioms can be turned into basic axioms.

\section{A game-theoretic interpretation of Tait calculus}\label{game}

To a sequent $\Gamma$ we associate a game $G(\Gamma)$ between Eloisa and Abelard.   Eloisa (resp.\ Abelard) is trying to show that $\Gamma$ is true (resp.\ false). Eloisa plays subformulas of $\E$-formulas. Abelard plays subformulas of $\A$-formulas. The game is designed so that the premise of an $\E$-inference 
\[ \ru{\Gamma, \vp_i}{\Gamma,\E \vp_i}{} \]
translates into a move $\vp_i$ by Eloisa. On the other hand,  each premise of an $\A$-inference 
\[ \ru{ \Gamma, \vp_{i} \ \ (i) }{\Gamma,\A\vp_i}{} \]
translates into a move $\vp_i$ by Abelard. Cut-free proofs of $\Gamma$  will then be interpreted as winning strategies for Eloisa in $G(\Gamma)$. Clearly, the rules of the game should not depend on a particular cut-free proof of $\Gamma$, if any. A simple solution is to regard the conclusion of an $\A$-inference  as a move by Eloisa. Essentially, Eloisa decides who makes the next move and which $\A$-formulas Abelard is allowed to challenge. 

\subsection{Games}

Although we are interested in a particular kind of game, it will be useful to consider a fairly general  game-theoretic framework. 

\begin{definition}
Let $X^{<\N}$  denote the set of finite sequences of elements of $X$. By $X^{\N}$  we denote the set of all (total) functions $f\colon \N\to X$. 
We denote sequences by $p,q,\ldots, \sigma,\tau, \ldots$ We use $\sigma\subseteq\tau$ for  initial segment, $\sigma\subset\tau$ for proper initial segment,  $\sigma\conc \langle x\rangle $ and $\sigma\conc\tau$ or simply $\sigma x$ and $\sigma\tau$ for concatenation.   
The empty sequence is also written $\lr$. We will often denote a finite sequence by $x_1\cdots x_n$.  The empty sequence is covered by the case $n=0$. Conversely,  we use  $x_0\cdots x_n$ to denote non-empty sequences. By writing $x_1x_2\cdots$ or $x_0x_1\cdots$ we denote finite and infinite sequences at once. The length of a finite sequence $\sigma$ is denoted by $|\sigma|$.   A (rooted) tree $T$ on $X$ is a non-empty subset of $X^{<\N}$ closed under initial segments. A leaf of $T$ is a node of $T$ with no proper extension in $T$. A branch node of $T$ is any node of $T$ which is not a leaf.
A tree $T\subseteq X^{<\N}$ is well-founded if for every $f\colon \N\to X$ there exists $n$ such that $f(0)\cdots f(n)\notin T$.  The height of a well-founded tree  $T$ is $\height(T)=\height(\lr)$, where $\height(\sigma)=\sup\{\height(\sigma x)+1\mid \sigma x\in T\}$ for every $\sigma\in T$. In general,  a (rooted) tree is a non-empty partial order $(T,\preceq)$ with a minimum element, the root, denoted $\lr$, such that $\{\sigma\in T\mid \sigma\preceq \tau\}$ is linearly ordered by $\preceq$ for every $\tau\in T$. A tree $T$ is well-founded if for every $f\colon \N\to T$ there exists $n$ such that $f(n)\not\prec f(n+1)$, where $\sigma\prec \tau$ is short for $\sigma\preceq\tau$ and $\sigma\neq\tau$. The height of a well-founded tree $T$ is $\height(T)=\height(\lr)$, where $\height(\sigma)=\sup\{\height(\tau)+1\mid \tau \in T\land \sigma\prec\tau\}$ for every $\sigma\in T$.
\end{definition}

\begin{definition}[games]
A \emph{game} is a tuple $G=(X,E,A,W)$, where 
\begin{itemize}[leftmargin=5mm]
	\item $X$ is a set (the set of moves),
	\item $E\cup A \subseteq X^{<\N}\times X$ (legal moves) with $\dom(E)\cap\dom(A)=\emptyset$ ($E$ for Eloisa and $A$ for Abelard), 
	\item and $W\subseteq P_G$ (the set of winning plays).
\end{itemize}
Here, $P_G\subseteq X^{<\N}$ (the set of plays) is the set of finite sequences $a_1\cdots a_n$ such that $(a_1\cdots a_i, a_{i+1})\in E\cup A$ for every $i<n$.
\end{definition}

\begin{definition}[strategy tree]
A {\em strategy tree} for Eloisa is a rooted tree $T\subseteq P_G$ such that for every $p = a_1\cdots a_n \in T$:
\begin{itemize}[leftmargin=5mm]
	\item if $p\in W$ then $p$ is a leaf;
	\item if $p\notin W$ and $p\in\dom(E)$, then there exists a unique $a$ such that $pa\in T$;
	\item if $p\notin W$ and $p\in\dom(A)$, then $p a\in T$ for every $a$ such that $(p,a)\in A$.
\end{itemize}
We say that $T$ is {\em non-losing} if $p\in W$ whenever $p\in T$ is a leaf. We say that $T$ is {\em winning} if $T$ is non-losing and well-founded.
\end{definition}

\begin{remark}
A strategy tree $T$ is non-losing iff for every $p\in T\setminus W$ we have $p\in\dom(E)\cup\dom(A)$.
\end{remark}

\begin{definition}
Let $f\colon X^{<\N}\to X$. Let  $T(f)$ be the tree  of all plays $a_1\cdots a_n\in P_G$ such that for every $i<n$ the initial segment $a_1\cdots a_i\notin W$ and $a_{i+1}=f(a_1\cdots a_{i})$  whenever $a_1\cdots a_i\in\dom(E)$. 
\end{definition}

\begin{definition}[strategy]\label{winning strategy}
(i) A  {\em strategy} for Eloisa  is a function $f\colon X^{<\N}\to X$ such that $T(f)$ is a  strategy tree.

(ii) A strategy for Eloisa $f$ is {\em non-losing} if $T(f)$ is non-losing. 

(iii) A  {\em winning strategy}  for Eloisa is a strategy $f$ such that  for every infinite sequence   $a_1a_2\cdots a_n\cdots$  there is an $n$ such that $a_1\cdots a_n\in P_G$ and  one of the following holds:
\begin{itemize}[leftmargin=5mm]
	\item $a_1\cdots a_n\in W$;
	\item $a_1\cdots a_n\in\dom(E)$ and $a_{n+1}\neq f(a_1\cdots a_n)$; 
	\item  $a_1\cdots a_n\in\dom(A)$ and $(a_1\cdots a_n,a_{n+1})\notin A$.
	\end{itemize} 
\end{definition}

\begin{definition}
The game $G$ is non-losing for Eloisa if $p\in P_G\setminus W$ implies $p\in\dom(E)\cup\dom(A)$. In  this case, every strategy tree  for Eloisa is non-losing.
\end{definition}

\begin{proposition}
Let $f$ be a strategy for Eloisa in the game $G$. Then $f$ is winning iff $T(f)$ is winning. If $G$ is non-losing for Eloisa, then $f$ is winning iff $T(f)$ is well-founded.  
\end{proposition}
\begin{proof}
Exercise. 
\end{proof}

\subsection{Games and ordinals}

We define the length of a winning strategy $f$ for Eloisa in $G(\Gamma)$ as the height of the  strategy tree $T(f)$. Let $\Ord$ denote the class of all ordinals.

\begin{definition}
	Let $T$ be a tree. A {\em height function}  on $T$ is a function $\hh\colon T\to \Ord$ such  $\hh(\tau)<\hh(\sigma)$ whenever $\sigma\prec \tau$ for all $\sigma, \tau\in T$.
\end{definition}

Notice that if $\hh$ is a height function on a tree $T$, then $T$ is well-founded and $\height(T)\leq \hh(\lr)$. 

\begin{proposition}
	Let $G$ be a game and $f$ be a strategy for Eloisa in $G$. Then $f$ is winning iff $T(f)$ is non-losing and there is a height function  on $T(f)$. If $G$ is non-losing for Eloisa, $f$ is winning iff $T(f)$ has a  height function.
\end{proposition}
\begin{proof}
	Obvious.
\end{proof}

\subsection{Tait  games} 

Let us fix a collection of sequents to be regarded as axioms. Throughout the paper, this is done by specifying a set $S$ of true literals so that a sequent $\Gamma$ is an axiom if $\Gamma\cap S\neq\emptyset$.  

\begin{definition}[Tait games]
Let $G(\Gamma)=(X,E,A,W)$ be the  game defined by:
\begin{itemize}[leftmargin=5mm]
\item $X$ consists of pairs $(\ell,\vp)$, where $\ell\in\{g,q,r\}$ and  $\vp$ is an $\A$-formula when $\ell=q$;
\item $E$ consists of pairs $(a_1\cdots a_n, a)$, where $n=0$ or $a_n$ is not of the form $(q,\vp)$ and $a$ is either $(g,\vp_j)$ with  $\E\vp_j\in \Gamma_n$ or  $(q, \vp)$ with $\vp\in \Gamma_n$;
\item $A$ consists of pairs  $(a_0\cdots a_n, a)$, where $a_n$ is $(q,\A\vp_j)$ and $a=(r,\vp_j)$ for some $j$;
\item $W=\{a_1\cdots a_n\in P \mid  \Gamma_n \text{ is an axiom}\}$,
 \end{itemize}
where $P$ is the set of plays and $\Gamma_n=\Gamma\cup\{\vp\mid (\ell,\vp)\in \{a_1,\cdots,a_n\} \}$.
\end{definition}

\begin{definition}
A move $(g,\vp)$ by Eloisa  is called a guess. We write $\guess \vp$. A move by Eloisa of the form $(q,\vp)$ is called a query. We write $\query \vp$.  Each move by Abelard is a reply. We write $\reply  \vp$. We use $\vp$ to denote a move of the form  $\guess \vp$ or $\reply  \vp$.  A move $\vp$ (a guess or a reply) is  {\em winning} if $\vp$ is a true literal.	
\end{definition}

\begin{remark}
Suppose $a_1\cdots a_n$ is a play in $G(\Gamma)$ and $a_n$ is a move by Eloisa of the form $\query \vp$. Then $\vp\in\Gamma$ or there is some $0<k<n$ such that $a_k$ is either $\guess \vp$ or $\reply  \vp$.
\end{remark}

\begin{remark}
A play in $G(\Gamma)$ looks like
\begin{multline*}  \cdots \quad \E\A\psi_{ij} \quad \cdots\quad  \A\vp_i \quad \cdots  \quad  \query \A\vp_i  \qquad   \reply \vp_2 \\ \guess \A\psi_{3j} \qquad   \query \A\psi_{3j} \qquad \reply \psi_{35} \qquad \guess \A\psi_{7j}     \quad \cdots 
\end{multline*}
Abelard can only play after a move by Eloisa of the form $\query \A \vp_i$ and try a counterexample $\vp_i$ to $\A \vp_i$.
\end{remark}

\begin{remark}\label{non-losing}
If $\Gamma$ contains either an $\E$-formula or an $\A$-formula, then $G(\Gamma)$ is non-losing for Eloisa.
\end{remark}

\subsection{Formulas as trees}

The following terminology will come in handy.

\begin{definition}
	A tree with {\em polarity} is a tree $T$ along with a subset $P\subseteq T$. Nodes $\sigma$ and $\tau$ have the same polarity if they are both in $P$ or $T\setminus P$. We usually omit the reference to $P$.
\end{definition}

\begin{definition}[depth]
	Let $T$ be a tree with polarity. The \emph{depth} of a node in $T$ is defined by:
	\begin{itemize}[leftmargin=5mm]
		\item $\dt(\lr)=0$;
		\item $\dt(\sigma i)= \dt(\sigma)$ if they have the same polarity, $=\dt(\sigma)+1$ otherwise;
	\end{itemize} 
	
A node $\sigma$ is {\em odd} (resp.\ {\em even}) if so is $\dt(\sigma)$. Say that $\sigma$ is \emph{minimal} if it  is a minimal node of depth $\dt(\sigma)$. A node $\tau$ is a \emph{minimal extension} of $\sigma$, written  $\sigma\subsetm \tau$, if $\sigma\subseteq \tau$, they are both minimal and $\dt(\tau)=\dt(\sigma)+1$. 
\end{definition}

\begin{remark}[polarity vs parity]
	Two nodes $\sigma$ and $\tau$ have the same polarity iff $\dt(\sigma)$ and $\dt(\tau)$ have the same parity. 
\end{remark}

\begin{definition}[formulas as trees]
We think of a formula $\vp$ as given by a  well-founded rooted tree $T$ with polarity $P$  along with an extra labeling $\ell$ of all leaves with literals. Leaves are required to be minimal:  if $\sigma i$ is a leaf,  $\dt(\sigma i)=\dt(\sigma)+1$. $\E$-formulas (resp.\ $\A$-formulas) have positive (resp.\ negative) polarity. We thus regard a branch node $\sigma$ in $P$ as having positive polarity. 

Let $\vp$ be a formula given by $T,P,\ell$. For every leaf $\sigma\in T$, the subformula $\vp_\sigma$ is the literal $\ell(\sigma)$. For every branch node $\sigma\in T$,  the subformula $\vp_\sigma$ is the positive formula $\E_{i\in I_\sigma} \vp_{\sigma i}$ if $\sigma\in P$, or the negative formula $\A_{i\in I_\sigma}\vp_{\sigma i}$ if $\sigma\notin P$, where $I_\sigma=\{i\mid \sigma i\in T\}$. In particular, $\vp$ is either a literal (in case $T$ consists only of the root $\lr$), a positive formula if $\lr\in P$, a negative formula otherwise. 
\end{definition}

\begin{convention}
We usually write $T(\vp)$ without reference to the polarity.  If we declare $\vp$ to be positive (negative), then it is understood that $P$  contains (does not contain) the root. 	
\end{convention}

\begin{remark}
Let	$\sigma\in T=T(\vp)$. Then $\neg\vp_\sigma= (\neg\vp)_\sigma$, where $T(\neg\vp)$ is the tree $T$ with polarity $T\setminus P$ and leaves $\sigma$ labeled by $\neg \ell(\sigma)$. In other words, every node $\sigma\in T$ identifies both the subformula  $\vp_\sigma$ of $\vp$ as well as the subformula $\neg\vp_\sigma$ of $\neg\vp$.
\end{remark}

\begin{remark}
In general, a leaf of $T(\vp)$ may be odd or even.
\end{remark}
\begin{remark}
	Suppose $\vp$ is an $\E$-formula and $\sigma\in T(\vp)$ is a branch node. Then 
	\begin{itemize}[leftmargin=5mm]
		\item  $\sigma$ is even iff $\vp_\sigma$  is an $\E$-formula, 
		\item  $\sigma$ is odd iff $\vp_\sigma$  is an $\A$-formula.
	\end{itemize}
\end{remark}

\subsection{Partial strategies for Abelard}

\begin{convention}[duality]
A duality principle applies to 	all definitions and properties marked with  (\**). The dual statements are obtained by  interchanging odd with even and  $\E$ with $\A$.  
\end{convention}

\begin{remark}[\**]
	Suppose $\vp$ is an $\E$-formula and let $p$ be a play in $G(\Gamma)$. Let $\tau\in T(\vp)$ be minimal. Suppose $\vp_\sigma$ is a move in $p$ with  $\sigma\subseteq\tau$. Then there is a unique pair of nodes $\zeta_0\subsetm\zeta_1$ with $\zeta_1 \subseteq\tau$ such that $\sigma\in (\zeta_0,\zeta_1]$. Moreover, if $\zeta_0$ is odd, then $\vp_\sigma$ is a reply by Abelard. If $\zeta_0$ is even, then $\vp_\sigma$ is a guess   by Eloisa. Similarly, suppose $\query \vp_\sigma$ is a move in $p$ with $\sigma \subset\tau$.  Then there is a unique pair of nodes $\zeta_0\subsetm\zeta_1$ with $\zeta_0$ odd and $\zeta_1 \subseteq\tau$ such that $\sigma\in [\zeta_0,\zeta_1)$. 
\end{remark}

\begin{definition}[\**]\label{odd}
	Let $T$ be a tree with polarity and $S\subseteq T$. Then $S$ is {\em odd}  if 
	\begin{enumerate}[label=(\roman*), leftmargin=0.8cm]
		\item every $\sigma\in S$ is minimal;
		\item if $\sigma \subsetm\tau$ and $\tau\in S$, then $\sigma\in S$;
		\item if $\sigma\in S$ is odd, then there is a unique $\tau\in S$ such that $\sigma\subsetm\tau$.
    \end{enumerate}	
\end{definition}
\begin{example}
	$S_0=\emptyset$ is even and $S_1=\{\lr\}$ is odd.	
\end{example}

\begin{definition}[\**]\label{partial strategy}
	Suppose $\vp$ is an $\E$-formula. A {\em partial strategy} for Abelard in $G(\vp)$ is a subset $S$ of $T(\vp)$ such that for every odd $\sigma\in T(\vp)$ there exists at most one $i$ such that $\sigma i\subseteq \tau$ for some $\tau\in S$. A play  $p$  in $G(\Gamma, \vp)$ {\em agrees} with $S$, or 
	Abelard  is said to {\em follow  $S$} in $p$, if  for every move $\reply \vp_\sigma$ in $p$ there is a $\tau\in S$ such that $\sigma\subseteq\tau$.
\end{definition}

\begin{lemma}[\**]\label{odd is strategy}
An odd subset of $T(\vp)$, where $\vp$ is an $\E$-formula, is a partial strategy for Abelard.
\end{lemma}

\begin{lemma}[\**]\label{odd strategy}
	Suppose $\vp$ is an $\E$-formula. Let $p$ be a play in $G(\Gamma,\vp)$ and suppose Abelard follows  $S$ in $p$, where $S$ is odd.  Then 
\begin{enumerate}[label=$(\arabic*)$, leftmargin=0.8cm]	
\item for every Abelard move $\reply \vp_\sigma$ there is a unique pair of nodes $\zeta_0, \zeta_1$ in $S$ such that  $\zeta_0$ is odd, $\zeta_0\subsetm\zeta_1$,  and $\sigma\in (\zeta_0,\zeta_1]$;
	
\item if $S$ is non-empty, for every Eloisa move of the form $\guess \vp_{\sigma}$ there is a unique $\tau\in S$ such that $\tau$ is even,  $\tau\subset \sigma$ and $\dt(\sigma)\leq \dt(\tau)+1$.
\end{enumerate}
\end{lemma}
\begin{proof}
(1) Suppose $\sigma$ has an extension in $S$. Note that $\sigma$ is not the root.  By Definition \ref{odd}\ (i) and (ii) we can find a pair of nodes $\zeta_0\subsetm \zeta_1$ in $S$ such that $\sigma\in(\zeta_0,\zeta_1]$. Clearly, $\zeta_0$ is unique. On the other hand,  if $\zeta_0$ is even, then  $\vp_\sigma$ cannot be a reply by Abelard. Then $\zeta_0$ is odd and so the pair is unique by Definition \ref{odd}\ (iii), as desired. 
	
(2) First, by Definition \ref{odd}\ (i) and (ii), if $S$ is non-empty then $\lr\in S$. Since $\vp_\sigma$ is a guess, there is a unique even minimal node $\tau\subset\sigma$  such that   $\dt(\sigma)=\dt(\tau)$ or $\dt(\sigma)=\dt(\tau)+1$.  If $\tau\neq\lr$, then $\vp_\tau$ must be an earlier reply by Abelard. By part (1), $\tau\in (\zeta_0,\zeta_1]$ for a unique pair $\zeta_0\subsetm \zeta_1$ in $S$. By minimality, we must have $\tau=\zeta_1$. Hence $\tau\in S$. 
\end{proof}

\section{Debates} \label{debates}

A strategy  for Eloisa in $G(\Gamma)$ can be obtained by playing a strategy $f$  for Eloisa in  $G(\Gamma,\vp)$ against  a strategy $g$ for Eloisa in the \emph{dual} game $G(\Gamma,\neg\vp)$.  The idea is to  use $f$ and $g$ to compute a \emph{debate}, namely, an  alternating sequence of plays in $G(\Gamma,\neg\vp)$ and $G(\Gamma,\vp)$ so that  all moves by Abelard in each play are determined by some earlier move by Eloisa in a previous play. 

To get an idea of how a debate unfolds, consider a formula $\vp=\E\vp_i$ and assume for simplicity that each  $\vp_i$ is an $\A$-formula or a literal.  We start by  playing in $G(\Gamma,\neg\vp)$.  We wait until  Eloisa  plays $\query \neg\vp$. Abelard does not know how to reply to this move as yet.  We then switch to $G(\Gamma,\vp)$ and wait until Eloisa plays   $\query \vp_i$ or a winning $\guess\vp_i$ form some $i$. (Notice that in the second case $\vp_i$ must be a literal.) We resume the game in $G(\Gamma,\neg\vp)$. Abelard can now answer with $\reply \neg\vp_i$. And so on and so forth. Meanwhile,  every move on the $\Gamma$ side, from either $G(\Gamma,\vp)$ or $G(\Gamma,\neg\vp)$, is a move in $G(\Gamma)$. 
The challenge is to show that if $f$ and $g$ are winning strategies,  either player will eventually play a true literal from $\Gamma$. 

Essentially, a debate is a sequence of plays $p_0p_1\cdots$, where $p_{2n}$ is a play in $G(\Gamma,\neg\vp)$ and $p_{2n+1}$ is a play in $G(\Gamma,\vp)$. The last move in every play $p_n$, if $p_{n+1}$ is defined, is either a winning guess   or a query. For $n$ odd, such move is of the form $\guess\vp_\sigma$ or $\query \vp_\sigma$. Dually for $n$ even.  Moreover, each play $p_{n+1}$, where $n>0$, extends some previous play $p_{i_{n}}$ with $i_{n}<n$ whose last move is a query.  If $n+1$ is odd and $p_{i_{n}}$ ends with $\query \vp_\sigma$, then $p_{n}$ will end with $\guess\neg\vp_\tau$ or $\query \neg\vp_\tau$  for some $\tau$ with $\sigma\subset\tau$. Thus in $p_{n+1}$ Abelard can answer to $\query \vp_\sigma$ with $\reply \vp_{\sigma i}$, where $\sigma i\subseteq\tau$. Dually for $n+1$ even.
We thus have a \emph{pointer} from $n$ to $i_{n}=m<n$ such that $p_{n+1}$ extends the play $p_m$ thanks to some move in the play $p_n$.   Coquand's termination argument crucially relies on the {\em geometry} of the pointer sequences associated to a debate, which he calls interaction sequences.

\begin{definition}[interaction sequences]
A {\em pointer sequence} is a sequence of natural numbers $i_0i_1i_2\cdots$ such that $i_0=0$ and $i_{n+1}<n+1$ for every $n$. 
An \emph{interaction sequence} is a pointer sequence  $i_0i_1i_2\cdots$ 
such that \[ i_{n+1}\in V(n+1)=\{n\}\cup V(i_n), \] where the sequence of finite sets  $V(n)$ is obtained by letting  $V(0)=\emptyset$.
\end{definition}
\begin{remark}
In any interaction sequence,	we must have $i_1=0$, $i_2=1$, $i_3\in\{0,2\}$, and so forth. Moreover, $i_{n+1}<n+1$ have opposite parity. 
\end{remark}

\begin{lemma}\label{int}
Let $i_0i_1i_2\cdots$ be an interaction sequence. The following holds:
	\begin{itemize}[leftmargin=5mm]
		\item if $m\in V(n)$, then $V(i_m)\subseteq V(n)$;
		\item if $k,m\in V(n)$ and $k<m$, then $k\in V(i_m)$.
\end{itemize}
\end{lemma}

\begin{definition}
	Given a pointer sequence $i_0i_1i_2\cdots$, define \[ W(n+1)=\{n,i_n\}\cup W(i_n),\]
	by letting  $W(0)=\emptyset$. Note that $W(n)=V(n)\cup\{i_m\mid m\in V(n)\}$.
\end{definition}

Let $\party(n)$ be the parity of $n$, i.e., $\party(n)$ is $1$ if $n$ is odd and $0$ otherwise.  
\begin{note}
In the context of  Hyland-Ong games \cite{HO00}), where plays are essentially conceived as abstract debates, the set $W(n)$ corresponds to the {\em view} of player $\party(n)$. The move $p_n$ is a {\em response} to  the  move $p_{i_n}$. Notice that the latter is a move by the opponent player $\party(i_n)$. Then $V(n)=W(n)\cap\{m\mid \party(m)\neq\party(n)\}$. In other words, $i_n\in V(n)$ ensures that $\party(n)$ does not respond to their own moves.  

\end{note}

The following  notion of  debate will be sufficient to carry out Coquand's proof of cut elimination (see Theorem \ref{cut strategy}).

\begin{definition}[debates]\label{debate}
Suppose $\vp$ is an $\E$-formula.	Let $f_0$ be a strategy for Eloisa in $G(\Gamma,\neg\vp)$ and $f_1$ be a strategy for Eloisa in $G(\Gamma,\vp)$.\footnote{Note that both strategies are non-losing (cf.\ Remark \ref{non-losing}).} Let $X$ be the set of moves  of $G(\Gamma)$.   

A {\em quasi-debate}  is a function $d$ mapping every $p=a_1a_2\cdots\in \dom(d)=X^{<\N}\cup X^\N$ to a (possibly infinite) sequence 
\[  (d_0,\pGamma_0)\cdots (d_s,\pGamma_s)\cdots,\] 
where each $d_s$ is a finite sequence of the form $$(p_0,\sigma_0,i_0)\cdots (p_{n-1},\sigma_{n-1}, i_{n-1})\,  p_{n}$$ such that the following holds for every $m\leq n$:
\begin{enumerate}[label=(\arabic*), leftmargin=0.8cm]
	\item  $p_{m}\in T(f_{\party(m)})$\,  (in particular, $p_m$ is a play in $G(\Gamma,\vp)$ if $m$ is odd  and in $G(\Gamma,\neg\vp)$ if $m$ is even);
	\item $\sigma_m\in T(\vp)$ for  $m< n$; 
	\item $i_0\cdots i_{n-1}$ is an interaction sequence such that $p_{i_m}\subset p_{m+1}$ and  $\dt(\sigma_{i_m})<\dt(\sigma_m)$ for $0<m<n$,
\end{enumerate}
and for each $s$:
\begin{enumerate}[label=(\arabic*), leftmargin=0.8cm,start=4]
	\item $\pGamma_s\subseteq p$ is a play  in $G(\Gamma)$;
	\item $d_s\subset d_{s+1}$ and $\pGamma_s\subseteq \pGamma_{s+1}$ whenever defined.
\end{enumerate}

The play $p_n$ is called the {\em leading} play in $d_s$.

We say that a quasi-debate $d$ is {\em monotone} if $q\subseteq p$ implies $d(q)\subseteq d(p)$.

A quasi-debate $d$ is considered {\em non-losing}  if, for every $s$, the following conditions are equivalent:
\begin{itemize}[leftmargin=5mm]
\item the sequence  $d(p)$ is finite with $|d(p)|=s+1$;
\item  one of the following {\em winning} requirements holds for $\dot p_s=a_1\cdots a_l$: 
\begin{itemize}[leftmargin=2mm]
	\item $\dot p_s\in W$;
	\item both $\dot p_s$ and $p_n$ belong to $\dom(E)$, $f_{\party(n)}(p_n)=a$ is a legal move in $G(\Gamma)$, but either $p=a_1\cdots a_l$ or $a_1\cdots a_{l+1}\subseteq p$ and $a\neq a_{l+1}$;
	\item $\dot p_s\in\dom(A)$, but either $p=a_1\cdots a_l$ or $a_1\cdots a_{l+1}\subseteq p$ and $a_{l+1}$ is not a legal move in $G(\Gamma)$.
\end{itemize}
\end{itemize}

A {\em debate} is a non-losing monotone quasi-debate. 
We  say that $d$ is a debate for $f_0, f_1, \Gamma, \vp$.  We also say that $d(p)$ is a debate for $p$.
\end{definition}

\begin{lemma}\label{continuity}
Every debate is continuous. If $d(p)$ is finite then there is a finite initial segment $q\subseteq p$ such that $d(q)=d(p)$.  Indeed, if $|d(p)|=s+1$ then $d(\dot p_s)=d(p)$.
\end{lemma}	
\begin{proof}
Let $q=\dot p_s$. By monotonicity,  $d(q)\subseteq d(p)$. Suppose for a contradiction that $d(q)\subset d(p)$. Say $|d(q)|=t+1$ with $t<s$. Note that $\dot q_{t}=\dot p_{t}$, and so $\dot p_{t}\subseteq q\subset p$. But then the debate for $p$ should have concluded by stage $t$. The last claim follows from $d$ being non-losing. 
\end{proof}

\begin{note}
This  notion of debate is quite vague about the role of the nodes $\sigma_m$'s. In our construction of the canonical debate (Definition \ref{canonical debate})  we will make sure that for every $m\leq n$ the set $$S_m=\{\sigma_k\mid k\in W(m)\}\subseteq T(\vp)$$ is a partial strategy for Abelard in $G(\vp)$ if $m$ is odd, respectively $G(\neg\vp)$ if $m$ is even, and Abelard plays according to $S_m$ in the play $p_m$ (cf.\ Definition \ref{partial strategy}).  
\end{note}

\section{Debates of finite depth}\label{finite depth debates}

In this section, we detail an   analysis of Coquand's termination argument for debates of {\em finite depth}. This is sufficient for dealing with Peano Arithmetic. In other words, we show  that   every debate of  finite depth gives rise to  a winning cut strategy  (see Theorem \ref{termination} and Theorem \ref{cut strategy}).

We thus assume that formulas have finite depth, namely, for every formula $\vp$ there exists a natural number $\nu$, the depth of $\vp$,  such that $\dt(\sigma)\leq \nu$ for every $\sigma\in T(\vp)$.

\begin{definition}
An interaction sequence $i=i_0i_1i_2\cdots $ has {\em depth} at most $\nu$ if the length of $n\ i_n\ i_{i_n} \cdots 0$ is $\leq \nu+1$ for every $n$.
\end{definition}

\begin{lemma}
Let $d$ be a debate for $f_0, f_1, \Gamma, \vp$. Then the depth of any interaction sequence $i=i_0i_1i_2\cdots $ associated  to $d$ is  at most $\nu$, where $\nu$ is the depth of $\vp$.
\end{lemma}
\begin{proof}
By Definition \ref{debate},  we have $\dt(\sigma_{i_m})< \dt(\sigma_m)$ for every $m>0$. 
\end{proof}

\begin{definition}
A  finite sequence $u=(\alpha_0,i_0)(\alpha_1,i_1)\cdots(\alpha_{n-1},i_{n-1})\, \alpha_n$ is an {\em ordinal interaction sequence} if 
\begin{itemize}[leftmargin=5mm]
\item  $i_0i_1i_2\cdots$ is an interaction sequence for $n>0$;
\item each $\alpha_m$ is an ordinal;
\item  $\alpha_{m+1}< \alpha_{i_m}$  for every $0<m<n$.
\end{itemize}
The sequence $u$ is of depth at most $\nu$  if so is $i_0i_1i_2\cdots$. Also, $u$ is below $\alpha$, written $u<\alpha$, if all ordinals are $< \alpha$.
\end{definition}

\begin{definition}
Let 
\[ 	u = (\alpha_0,i_0)(\alpha_1,i_1)\cdots(\alpha_{n-1},i_{n-1})\, \alpha_n \text{\quad and \ }
v = (\beta_0,j_0)(\beta_1,j_1)\cdots\cdots (\beta_{m-1},j_{m-1})\, \beta_m \]
 be ordinal interaction sequences. Write $u \preceq v$ if
\begin{itemize}[leftmargin=5mm]
	\item  $n\leq m$,
	\item $\alpha_k=\beta_k$ and $i_k=j_k$ for every $k<n$,
	\item  $\beta_n\leq \alpha_n$.
\end{itemize}
\end{definition}

For the purposes of Theorem \ref{termination}, we wish to show  that ordinal interaction sequences of assigned depth are well-founded under extension, namely, there cannot be an infinite sequence $u_0 \prec u_1\prec \cdots \prec u_s\prec\cdots$ of ordinal interaction sequences of depth at most $\nu$, for any given $\nu$.\footnote{On the other hand, there are  infinite interaction sequences of bounded depth.   For example, let 
\[  0012325272\cdots \]
be the infinite sequence defined by letting $i_0=i_1=0$, $i_{2n+1}=2$ and $i_{2n+2}=2n+1$. One can easily see that this is an interaction  sequence of depth $\leq 4$.} This will require a few definitions and lemmas.

\begin{definition}
Let $i=i_0i_1i_2\cdots$ be a pointer sequence.  We say that $I \subseteq \dom(i)=\{0,1,\ldots\}$ is {\em isolated}  if $i_n\in I$ implies $n\in I$ for all $n$. We define $i-I$ to be the empty sequence if $I=\{0,1,\ldots\}$ or the unique pointer sequence $j$ such that  $i_{e(n)}=e(j_n)$, where $e(0)<e(1)<\cdots$ enumerates  $\{0,1,\ldots\}\setminus I$. 
\end{definition}

\begin{definition}
Let $i_0i_1i_2\cdots \cdots i_m\cdots $ be an pointer sequence. An {\em interval}  of $i$ is one of the form  $[i_m,m]$.   We say that $m$ has depth $\nu$ if the length of $m\ i_m\ i_{i_m} \cdots 0$ equals $\nu+1$. 
\end{definition}

What makes  interaction sequences special, among all pointer sequences, are the following three lemmas. To facilitate the verification of these lemmas,  it will be convenient to write $in$ for $i_n$ and let $h$ be such that $h0=0$ and $h(n+1)=n$. It then follows from the definition of interaction sequence  that for all $n$ the value of $in$ is of the  form $(hi)^s\, hn$ for some $s\in\N$, where in general by $f^s$ we denote the $s$-th iterate of $f$.

\begin{lemma}%[cf.\ {\cite[Lemma 1]{C95}}]
Let $i=i_0i_1i_2$ be an interaction sequence.  If $m\neq i_n$ for every $n>m$, then the interval $[i_m,m]$   is isolated.  In particular, if $i_0i_1i_2\cdots$ has depth at most $\nu$ and $m$ has depth $\nu$, then $[i_m,m]$ is isolated.
\end{lemma}
\begin{proof}
For the first part,  if $in>m$, there is nothing to show. Suppose now that $n>m$ and  $in\leq m$. Since $m$ cannot be the value of $in$, we must have $in<m$. We wish to prove that $in<im$.  By induction, we may safely assume that this holds for every $k$ with $m<k<n$ and $ik<m$.  Since $n>m$, there must be in the history of $in$, namely the sequence
\[ n> n_0=hn > \cdots > n_t=(hi)^t\, hn>\cdots > n_s=(hi)^s\, hn= in,\]
a  stage $t< s$ such that  $m\leq n_t$ and $n_{t+1}<m$. If $m=n_t$, then $n_{t+1}=him<im$, and so $in<im$. On the other hand, if $m<n_t$ then $in_t<m$. In fact, $hin_t=n_{t+1}<m$ and hence $in_t\leq m$. Again,  $m$ cannot be the value of $in_t$. We can now apply the induction hypothesis to $k=n_t$ and obtain $in_t<im$. Whence, as before, we get $in<im$.

The second part is obvious.
\end{proof}

\begin{lemma}[cf.\ {\cite[Lemma 1 and 2]{C95}}]
Let $i=i_0i_1i_2$  be an interaction sequence and  $I$ be a  union of isolated intervals. Then $I$ is isolated  and $i-I$ is either empty or an interaction sequence. 
\end{lemma}
\begin{proof}
Clearly, $I$ is isolated. Suppose $i-I$ is not empty.   Enumerate $E=\{0,1,\ldots\}\setminus I$ by $e(0)<e(1)<\cdots$.    Let us write $jn$ for $j_n$.   By definition, 
\begin{itemize}[leftmargin=5mm]
\item $ie(n)=e(jn)$.
\end{itemize}
We want to show that $j=i-I$ is also  an interaction sequence. We will prove that if 
 $ie(n)$ is of the form $e(n)_s=(hi)^s\, he(n)$ for some $s$, then  $jn$ is of the form $n_t=(hj)^t\, hn$ for some $t\leq s$.

Note that isolated intervals  are either disjoint or strictly included one into the other. Then we can partition $\{0,1,\ldots\}$ into non-empty intervals $E_0\cup I_0\cup E_1\cup \cdots$, where $E_0<I_0<E_1<\cdots$, each $I_k$ is the union of adjacent isolated intervals and $I= I_0\cup I_1\cup \cdots$.  Here, $J<K$ means that $n<m$ for all $n\in J$ and $m\in K$. We refer to the intervals $E_k$ as to the (connected) components of $E$.

If $e(n)$ and $ie(n)=e(n)_s$ are in the same component, it is immediate to see by induction on $s$ that $jn=n_s$.  In this case $ie(n)$ and  $jn$ share the same history. 

Suppose that $e(n)$ and $ie(n)=e(n)_s$ are not in the same component. Suppose first that  $I=I_0$ and so $\{0,1,\ldots\}=E_0\cup I\cup E_1$ with  $ie(n)\in E_0$ and $e(n)\in E_1$. Say $I=[k,m]$. Since $e(n)>m$ and $ie(n)<k$, there is  a $t<s$ such that $e(n)_t\geq m$ and $e(n)_{t+1}<m$.  First case: $e(n)_t>m$. By induction on $t$ one can verify that $e(n)_t=e(n_t)$.  From $e(n)_{t+1}<m$ it follows that $ie(n)_t<m$. But $I$ is isolated, hence $ie(n)_t<k$.  In this case,   $e(n)_{t+1}=hie(n_t)=e(hjn_t)=e(n_{t+1})$. Hence, $e(n)_s=e(n_s)$ and $jn=n_s$. Again, $ie(n)$ and  $jn$ share the same history.  Second case: $e(n)_t=m$.  By induction on $t$ one can verify that $e(n)_t+1=e((jh)^tn)$. It is easy to check that $e(n)_{t+r}=(hi)^rm=k-1$, where $r>0$ is the number of isolated intervals in $I$. Then $e(n)_{t+r}=e(h(jh)^tn)=e((hj)^thn)=e(n_t)$. Therefore, $e(n)_s=(hi)^{s-t-r}e(n_t)=e(n_{s-r})$ and $jn=n_{s-r}$. In this case $jn$ has a shorter history than $ie(n)$. The general case, where $I$ is not an interval, is obtained by iterating the above argument.  
\end{proof}

\begin{lemma}
Let $u=(\alpha_0,i_0)(\alpha_1,i_1)\cdots  (\alpha_{n-1},i_{n-1})\, \alpha_n$ be an ordinal interaction sequence. If $I$ is a union if isolated intervals, then 
\[ u-I=(\alpha_{e(0)},j_0)\cdots \cdots (\alpha_{e(m-1)},j_{m-1})\, \alpha_n,\]
where $j=i-I$, is an ordinal interaction sequence. 
\end{lemma}
\begin{proof}
By the previous lemma, $j$ is empty or an interaction sequence. If $j$ is empty, then $u-I=\alpha_n$, which is an ordinal interaction sequence. Suppose we are in the second case. We only need to show  that $\beta_{k+1}<\beta_{j_k}$ for all $0<k<m$, where $\beta_k=\alpha_{e(k)}$ for $k<m$ and $\beta_m=\alpha_n$. Once again,   write $ik$ for $j_k$.  Consider the case $k+1<m$. Then
\[  \tag{$*$} \beta_{k+1}=\alpha_{e(k+1)}<\alpha_{i(e(k+1)-1)}. \]

First case:  $e(k+1)-1=e(k)$, and so the  ordinal on the right in $(*)$  is $\alpha_{ie(k)}$. By definition of $j$, $ie(k)=e(jk)$, and so $\alpha_{ie(k)}=\alpha_{e(jk)}=\beta_{jk}$, as desired. 

Second case:  $e(k+1)-1$ is the right end-point of an isolated interval in the union. In this case,  $e(k)+1= (ih)^s\, ihe(k+1)$ for some $s\in \N$. By induction on $s$ one can verify that $\alpha_{i(e(k+1)-1)}=\alpha_{ihe(k+1)}\leq \alpha_{(ih)^s\, ihe(k+1)}$ and  therefore the ordinal on the right in $(*)$ is $\leq\alpha_{e(k)+1}$. Now,  $\alpha_{e(k)+1}<\alpha_{ie(k)}$. As in the first case, $\alpha_{ie(k)}=\beta_{jk}$, and we are done. 

For $k+1=m$, that is $\beta_{k+1}=\alpha_n$, we have
\[ \beta_{k+1}<\alpha_{i(n-1)}. \]
There are two cases, according to whether $n-1$ belongs or not to $I$.  If $n-1\in I$, then   $[i(n-1),n-1]$ is among the isolated intervals.  Then $e(m-1)+1=(ih)^s\, ihn$ for some $s$. As before, one can verify that  $\alpha_{i(n-1)}=\alpha_{ihn}\leq \alpha_{e(m-1)+1}$. On the other hand,  $\alpha_{e(m-1)+1}<\alpha_{ie(m-1)}=\alpha_{e(j(m-1))}=\beta_{j(m-1)}=\beta_{jk}$, as required. Suppose now that $n-1\in E$. In this case, we must have $n-1=e(m-1)$. The conclusion follows as in the previous case. The proof is complete.
\end{proof}

\begin{theorem}[cf.\ {\cite[Proposition 3]{C95}}]\label{height}
Let $\nu$ be a natural number and $\alpha$ be an ordinal. The class tree $T(\nu)$ of ordinal interaction sequences of depth at most $\nu$ is well-founded. The height of the tree $T(\nu,\alpha)$ of ordinal interaction sequences $u<\alpha$   of depth at most $ \nu$ is at most $c_\nu(\alpha)$, where
 \begin{align*}
	c_0(\alpha) & = \alpha \cdot 2\\ 
	c_{\nu+1}(\alpha) &= 3^{\, c_\nu(\alpha)} 
\end{align*}  
\end{theorem}

\begin{proof}
The first claim follows from the second. Simply note that if $u=(\alpha_0,i_0) \alpha_1$ is an ordinal interaction sequence, then every extension of $u$ is below $\max(\alpha_0,\alpha_1)+1$. 

The second claim is by induction on $\nu$. Note that $T=T(\nu,\alpha)=\{\lr\}\cup\{ u\in T(\nu)\mid u \text{ is below } \alpha\}$  and  the height of $T$ is $\height(\lr)$.

Case $\nu=0$. The only interaction sequence of depth $\leq 0$ is $i_0$. Then $u\in T(0,\alpha)$ is of the form $\lr$ or $\alpha_0$ or $(\alpha_0,i_0) \alpha_1$. The height of $T(0,\alpha)$ is $\leq \alpha\cdot 2$.

%Case $N=1$. The only interaction sequences of depth $\leq 1$ are $i_0$ and $i_0i_1$.  The height of $T(1,\alpha)$ is $\leq 2\alpha\cdot 3=c_2(\alpha)$.\prec  

Case $\nu+1$. Suppose $T_\nu=T(\nu,\alpha)$ is well-founded of height $\rho$. We show that the height of $T_{\nu+1}=T(\nu+1,\alpha)$ is at most $3^{\rho}$. 
First, we map sequences $u'\in T_{\nu+1}$ to  sequences $u\in T_\nu$  by removing all isolated intervals $[i_m,m]$ where $m$ has depth $\nu+1$. That this operation has the desired effect is guaranteed by the foregoing three lemmas.  Second, to $u\in T_\nu$ we assign an ordinal $c(u)$ as follows.   Let $c(\lr)=3^{\rho}$. Given $u=(\alpha_0,i_0)(\alpha_1,i_1)\cdots(\alpha_{n-1},i_{n-1})\, \alpha_n$ in $T_\nu$, let
\[ c(u)=3 ^{\rho_0} + \cdots + 3 ^{\rho_{n-1}}+3^{\rho_{n}}\cdot 2, \]
where  $\rho_{m}$ is the height of $(\alpha_0,i_0)(\alpha_1,i_1)\cdots(\alpha_{m-1},i_{m-1})\, \alpha_m$ in $T_\nu$. Note that $c(u)$ is in Cantor normal form. It is not difficult to see that  $u\prec v$ in $T_\nu$ implies $c(u)>c(v)$. Now, we claim that if $u'\prec v'$ in $T_{\nu+1}$ then $c(u)>c(v)$. In fact, if $u=(\alpha_0,i_0)(\alpha_1,i_1)\cdots(\alpha_{n-1},i_{n-1})\, \alpha_n$, then either $u\prec v$ in $T_\nu$ or for some $m<n$ there exists $\beta_m$ with $\beta_m< \alpha_m$  such that 
\[  w= (\alpha_0,i_0)\cdots (\alpha_{m-1},i_{m-1})\, \beta_m \preceq v. \]
In the latter  case,  $c(u)>c(v)$ as well.   We have just shown that if $T_\nu$ has height $\leq c_\nu(\alpha)$, then $u'\mapsto c(u)$ is a height function on  $T_{\nu+1}$ witnessing  $\height(T_{\nu+1})\leq c_{\nu+1}(\alpha)$.
\end{proof}

\begin{theorem}[termination]\label{termination}
Suppose $\vp$ is an $\E$-formula of finite depth.	Let  $f_0$ be a \wse in $G(\Gamma,\neg\vp)$ and $f_1$ be a \wse in $G(\Gamma,\vp)$. Let $p=a_1a_2\cdots $ be a sequence of moves in $G(\Gamma)$ and $d(p)$ be a debate for $p$. Then $d_s$ is undefined for some $s$.  
\end{theorem}

\begin{proof}
Suppose  $d_s$ is defined for every $s$ and let $\nu$ be a bound on the depth of $\vp$. We are going to define a descending sequence in $T(\nu)$, for the desired contradiction. Indeed, let $\hh_i$ be a height function on $T(f_i)$ and $\alpha=\max(\hh_0(\lr),\hh_1(\lr))$. We construct a descending sequence $(u_s)_s$ in $T(\nu,\alpha+1)$.  At stage $s$, if $d_s=e_s\, p_n=(p_0,\sigma_0, i_0)\cdots (p_{n-1},\sigma_{n-1}, i_{n-1})\,  p_{n}$,  we let  $u_s$ be $(\alpha_0,i_0)(\alpha_1,i_1)\cdots(\alpha_{n-1},i_{n-1})\, \alpha_n$, where $\alpha_m=\hh_{\party(m)}(p_m)$. By Definition \ref{debate}, each $p_m\in T(f_{\party(m)})$ and  $p_{i_m}\subset p_{m+1}$. It thus follows that   $\alpha_{m+1}<\alpha_{i_m}$ for every $m<n$, and so $u_s$ is indeed an ordinal interaction sequence. 
We now check that  $u_{s}\prec u_{s+1}$. In fact, $u\subset u_{s+1}$. But this immediately follows from $d_s\subset d_{s+1}$. In fact, either $d_{s+1}=e_s\, q_n$ with $p_n\subset q_n$, or $d_{s+1}=e_s\, (p_n,\sigma_n,i_n)\, p_{n+1}$ with $p_{i_n}\subset p_{n+1}$. 
\smallskip

\end{proof}

It is now clear how to define the cut strategy.
\begin{definition}[cut strategy]
Let $d$ be a debate for $f_0,f_1,\Gamma, \vp$, where $\vp$ is an $\E$-formula. The {\em cut strategy } $f=\cut(f_0,f_1)$ for Eloisa in $G(\Gamma)$ is defined so.  Given  $p\in X^{<\N}$, compute the final round $d_s$ in the debate for $p$, if any.  
Let $f(p)=f_{\party(n)}(p_n)$, where $p_n$ is the leading play in $d_s$. 
\end{definition}

\begin{theorem}[winning cut strategy]\label{cut strategy}
Let $d$ be a debate for $f_0,f_1,\Gamma, \vp$, where $\vp$ is an $\E$-formula of finite depth. 	If $f_0$ is a \wse  in $G(\Gamma,\neg\vp)$ and $f_1$ is a \wse in $G(\Gamma,\vp)$, then $f=\cut(f_0,f_1)$ is a winning strategy for Eloisa in $G(\Gamma)$. 
\end{theorem}
\begin{proof}
Let $p=a_1a_2\cdots$ be an infinite sequence of moves. By Theorem \ref{termination},  the debate for $p$ stops at some stage $s$.  Let $\pGamma_s=a_1\cdots a_l$.  Compute $f(a_1\cdots a_l)$. Note that $d_s$ is also the final round in the debate for $a_1\cdots a_l$ by Lemma \ref{continuity}.  By the winning requirements  and the definition of  cut strategy, either $a_1\cdots a_l\in W$, or $a_{l+1}\neq f(a_1\cdots a_l)$ where $a_1\cdots a_l\in\dom(E)$ and $f(a_1\cdots a_l)$ is a legal move in $G(\Gamma)$, or $a_1\cdots a_l\in\dom(A)$ and $a_{l+1}$ is an illegal move by Abelard. In other words, $f$ is winning (cf.\ Definition \ref{winning strategy}).
\end{proof}

\begin{problem}
Show that the cut strategy is winning in general. That is to say, provide a termination argument for debates involving formulas $\vp$ of infinite rank (a countable ordinal). It would suffice to show that there are no infinite ordinal interaction sequences $(u_0,i_0)(u_1,i_1)\cdots(u_{n},i_{n})\cdots$ such that $i_0i_1\cdots$ is bounded in the sense that  for every  strictly increasing $f\colon\N\to\N$  there exists $n$ such that $f(n)\neq i_{f(n+1)}$.
\end{problem}

\section{Computing debates}\label{computing}

We introduce our working notion of debate. The proof of the main theorem  (Theorem \ref{cut}) relies  on such  construction.

\begin{definition}[canonical debate]\label{canonical debate}	Suppose $\vp$ is an $\E$-formula.	Let $f_0$ be a strategy for Eloisa in $G(\Gamma,\neg\vp)$ and $f_1$ be a strategy for Eloisa in $G(\Gamma,\vp)$. Let $p=a_1a_2\cdots $ be a sequence of moves in $G(\Gamma)$.  We compute a debate  $d(p)$ for $p$ as follows. \smallskip
	
Construction:\medskip
	
Stage $s=0$. Let $d_s$ consist of   $p_0=\lr$ and $\pGamma_s=\lr$. \smallskip
	
	Stage $s+1$. Suppose we have already defined
	$d_s= \overbrace{(p_0,\sigma_0, i_0)\cdots (p_{n-1},\sigma_{n-1}, i_{n-1})}^{e_s}\,  p_n$
	and  $\pGamma_s=a_1\cdots a_l$.  \smallskip

	Case A. The play $\pGamma_s$ is winning. \texttt{STOP}. 
	
	Case B. The play $p_{n}$  is winning. 
	
	Case C.  It is Eloisa's turn to move in $p_{n}$. Let $f_{\party(n)}(p_{n})=a$. 
	
	Case D. It is Abelard's turn to move in $p_{n}$.\medskip

	Case $n>0$ with $n$ odd (the even case is dual).  \smallskip

	Case B.1. The  last move is a winning  $\guess \vp_\sigma$ by Eloisa and there is a unique $m\in V(n)$ such that $\sigma_m\subsetm\sigma$. Let $\sigma_m j\subseteq\sigma$. Set
	\begin{align*}
		d_{s+1} & = e_s\, (p_{n},\sigma,m)\, p_{n+1}, & 
		p_{n+1} & = p_m\neg\vp_{\sigma_mj}, & 
		\pGamma_{s+1} & = \pGamma_s.
	\end{align*}
	
	Case B.2. Else, \texttt{EXIT}.

	%Case C.1. The move $a$ is  an illegal move in $p_n$.  STOP. \smallskip
	
	Case C.1. The move $a$ is a legal move on the $\Gamma$ side. Recall $\pGamma_s=a_1\cdots a_l$.
	
	Case C.1.1. $p=a_1\cdots a_l$ or else $a_1\cdots a_{l+1}\subseteq p$ and $a\neq a_{l+1}$. \texttt{STOP}.
	
	Case C.1.2.  Else, let 
	\begin{align*}
		d_{s+1} & =  e_s\, (p_{n}a), &
		\pGamma_{s+1} & = \pGamma_s\, a. 
	\end{align*}
	
	Case C.2. $a$ is $\guess \vp_\sigma$ or $\query \vp_\sigma$. Let 
	\begin{align*}
		d_{s+1} & = e_s\, (p_{n}a), &
		\pGamma_{s+1} &= \pGamma_s. 
	\end{align*}
	
	Case D.1. The last move in $p_{n}$ is of the form $\query \gamma_\sigma$ for some $\gamma\in\Gamma$.  
	
	Case D.1.1. If $p=a_1\cdots a_l$ or $a_{l+1}$ is not a legal move, \texttt{STOP}. 
	
	Case D.1.2. Else, let
	\begin{align*}
		d_{s+1} &= e_s\, (p_{n}a_{l+1}), &
		\pGamma_{s+1} & = \pGamma_s\, a_{l+1}.
	\end{align*}

	Case D.2. The last move in $p_{n}$ is  $\query \vp_\sigma$ for some $\sigma$.   
	
	Case D.2.1. %Abelard can move in $p_{n}$, namely, there is $m\in W(n)$ such that $\sigma\subset\sigma_m$. 
	There is a unique $m\in V(n)$  with $m>0$ such that  $\sigma_{i_m}\subseteq \sigma \subset \sigma_m$. Let $\sigma j\subseteq\sigma_m$. Let 
	\begin{align*}
		d_{s+1} &= e_s\, (p_{n}\vp_{\sigma j}), &
		\pGamma_{s+1} & = \pGamma_s. 
	\end{align*}
	
	Case D.2.2. %Abelard cannot move in $p_{n}$.  
	There is a unique  $m\in V(n)$  such that $\sigma_m\subsetm \sigma$. Let $\sigma_m j\subseteq\sigma$.  Let
	\begin{align*}
		d_{s+1} & = e_s\, (p_{n},\sigma,m)\, p_{n+1}, &
		p_{n+1} & = p_m \neg\vp_{\sigma_m j}, &
		\pGamma_{s+1} & = \pGamma_s.
	\end{align*}
	
	Case D.2.3. Else, \texttt{EXIT}.
	\smallskip
	
	Case $n=0$. This case differs from $n>0$ in that case B  never occurs and case D.2 runs as follows.   \smallskip
	
	Case D.2. The last move in $p_{0}$ is  $\query \neg \vp$. 
	Let
	\begin{align*}
		d_{s+1} & = (p_{0},\lr,0)\, p_{1} & p_1 & =\lr, &
		\pGamma_{s+1} & =  \pGamma_s.
	\end{align*}  
	
	\smallskip
	
	End of construction.	
\end{definition}

\begin{remark}
	Suppose $d_s\downarrow$. Then $\pGamma_s$ is obtained from $q_0q_1\cdots q_n$  by deleting all moves of the form $(\ell,(\neg)\vp_\sigma)$, where $q_0=p_0$, $q_1=p_1$ and $q_{n+2}$ is unique such that $p_{n+2}=p_{i_{n+1}}q_{n+2}$. 
\end{remark}

\begin{lemma}[verification]\label{verification} The canonical debate is a debate.  
\end{lemma}

\begin{proof}
	By induction on $s$ we prove that, whenever $d_s$ is defined, we do not exit the construction and the following hold for every $m\leq n$: 
\begin{enumerate}[label=(\alph*), leftmargin=0.8cm] 
        \item  $p_{m}\in T(f_{\party(m)})$;
		\item  $\sigma_0=0$; for $m<n$ odd, $\sigma_m\in T(\vp)$ is an odd node and the last move in $p_m$ is either a winning  $\guess \vp_{\sigma_m}$ or   $\query \vp_{\sigma_{m}}$   (resp.\  for $m<n$ even, $\sigma_m$ is even and the last move in $p_m$ is either a winning  $\guess \neg\vp_{\sigma_{m}}$ or  $\query \neg \vp_{\sigma_{m}}$); 
		\item  $i_0\cdots i_{n-1}$ is an interaction sequence such that $p_{i_m}\subset p_{m+1}$ and $\sigma_{i_m} \subsetm\sigma_m$   for  $0<m<n$;
		\item $\pGamma_s\subseteq p$ is a play  in $G(\Gamma)$;
		\item  $d_s\subset d_{s+1}$ and $\pGamma_s\subseteq \pGamma_{s+1}$ whenever defined;
		\item if the last move in $p_n$ is on the $\Gamma$ side, then  $\pGamma_s$ ends with the very same move; if the last move in $\pGamma_s$ is a query, then $p_n$ ends with the very same move; 
		\item  $S_{m}=\{\sigma_k\mid k\in W(m)\}$ is odd for  $m\leq n $  odd (resp.\ even for  $m\leq n$  even);
		\item Abelard follows the partial strategy $S_{m}$ in $p_{m}$   and for $m<n$ he has no moves in $p_m$ if the last move is a query.	
	\end{enumerate}
	
Once this is proved, it easily follows that $d$ is a debate. In fact, (a)--(e) clearly imply (1)--(5). By construction, $d$ is monotone. On the other hand, by design,  we only stop when a winning requirement is met.  The verification is simple but tedious. 	We only treat the main cases, namely,  B and D.2. \smallskip
	
The cases $s=0$ and $n=0$ are straightforward. Suppose 
$$d_s= (p_0,\sigma_0, i_0)\cdots (p_{n-1},\sigma_{n-1}, i_{n-1})\,  p_{n},$$ 
and $n>0$ is odd (the even case is dual). By the induction hypothesis,  $S_{n}$ is odd and Abelard follows $S_{n}$ in $p_{n}$.  Also, by (b), (c) and (e), $\sigma_m\in S_{n}$ is  even  iff  $m\in W(n)$  is even iff $m\in V(n)$. Also note that $\lr\in S_n$ for $n>0$.\smallskip
	
Case B. We first show that we do not exit, namely, 
	\begin{itemize}[leftmargin=5mm]
		\item the last move is a winning guess by Eloisa of the form   $\guess \vp_\sigma$,
		\item there is a unique $m\in V(n)$ such that $\sigma_m\subsetm\sigma$.
	\end{itemize}
Since $\dot p_s$ is not winning, by (f) the last move in $p_n$ must be a winning guess  or a winning reply of the form $\vp_\sigma$. Suppose for a contradiction that  the last move in $p_n$ is by Abelard.  It follows  by Lemma \ref{odd strategy}\ (1), since $p_n$ agrees with  $S_n$, that there must be $m\in V(n)$ such that $\sigma=\sigma_m$. Note that $m>0$. By (b), $\neg\vp_{\sigma}$ is the last move by Eloisa in $p_m$, and so $\vp_\sigma$ is false, a contradiction. Hence the last move in $p_{n}$ is indeed by Eloisa. Therefore, by Lemma \ref{odd strategy}\ (2), there is   $m\in V(n)$ such that $\sigma_m\subset \sigma$ and $\dt(\sigma)\leq \dt(\sigma_m)+1$. Then  $\sigma_m\subsetm\sigma$ since $\sigma$ is a leaf. We claim that $m$ is unique. Suppose that for a different $k\in V(n)$ we have $\sigma_k\subsetm\sigma$. We can assume $k<m$. First, we must have $\sigma_k=\sigma_m$ since there is unique $\zeta$ such that $\zeta\subsetm\sigma$. By Lemma \ref{int}, $k\in V(i_m)$. Also, as $m>0$, by (c), we have $\sigma_{i_m}\subsetm\sigma_m$. It then follows that $k>0$ and $\sigma_{i_k}=\sigma_{i_m}$. We now have a contradiction because $\query \vp_{\sigma_{i_m}}$ is the last move in $p_{i_m}$ and yet Abelard can make a move since  $\sigma_k\in S_{i_m}$, against the induction hypothesis (h). % so much for that
	
We now verify (a)--(h). By construction, $\sigma_{n}=\sigma$ and $i_{n}=m$. 
	
\begin{enumerate}[label=(\alph*), leftmargin=0.8cm]	
\item $p_{n+1}$ is a play in $G(\Gamma,\neg\vp)$. In fact, the last move in $p_m$ is $\query \neg \vp_{\sigma_{m}}$ and $p_{n+1}$ extends $p_m$ in $G(\Gamma,\neg\vp)$ with a legal move by Abelard;
\item $\sigma_{n}=\sigma$ is odd. In fact, $\sigma_m$ is even and $\dt(\sigma)=\dt(\sigma_m)+1$. The last move in $p_{n}$ is the true literal $\vp_{\sigma_{n}}$;
\item $i_{n}=m\in V(n)$, and so $i_0i_1\cdots i_n$ is an interaction sequence. By construction, $p_{i_{n}}\subset p_{n+1}$ and  $\sigma_{i_{n}}\subsetm\sigma_{n}$;  
\item $\pGamma_{s+1}=\pGamma_s$;
\item $d_s\subset d_{s+1}$ by (a) and  $\pGamma_s\subseteq\pGamma_{s+1}$ by (d);
\item clear;
\item $S_{n+1}$ is even. In fact, $S_{n+1}=\{\sigma,\sigma_m\}\cup S_m$, and $S_m$ is even by induction. Clearly, Definition \ref{odd}\ (i) holds. For Definition \ref{odd}\ (iii), since $\sigma_m\subsetm\sigma$,  we only need verify that $\sigma$ is unique in $S_{n+1}$ with this property.  If $m=0$, then $S_m=\emptyset$, and so this is clearly the case. If $m>0$ and there is a $\tau\in S_m$ such that $\sigma_m\subsetm\tau$, then $\query \vp_{\sigma_m}$ cannot be the last move in $p_m$ by (i). Finally, if $\tau \subsetm\sigma_m$ then $m>0$ and $\tau=\sigma_{i_m}\in S_m$, and thereby Definition \ref{odd}\ (ii) holds;
\item it is clear that Abelard follows $S_{n+1}$ in $p_{n+1}$. On the other hand, the last move in $p_n$ is a literal and so there is nothing else to verify.
	\end{enumerate}
	
	Case D.2.  Let us show that we do not exit.  	
	
	Case D.2.1. Suppose Abelard can make a move in $p_n$ according to the partial strategy $S_n$, that is, $\sigma\subset\tau$ for some $\tau\in S_n$.  We then claim that  there is a unique $m\in V(n)$ such that $\sigma_{i_m}\subseteq \sigma \subset \sigma_m$.  By the closure properties of odd sets (Definition \ref{odd}),  there is a pair $\zeta_0\subsetm\zeta_1$ in $S_n$ such that $\sigma\in [\zeta_0,\zeta_1)$. Then $\sigma_0$ is odd and the pair is unique. Also, there is an $m\in V(n)$ such that $\zeta_1=\sigma_m$ and $\zeta_0=\sigma_{i_m}$. The proof that $m$ is unique is as before.  
	
	Case D.2.2. Suppose on the contrary that Abelard cannot make a move in $p_n$ according to $S_n$. We show that  for some unique  $m\in V(n)$ we have $\sigma_m\subsetm\sigma$.   We first claim that $\sigma$ is minimal. If not,  $\vp_\sigma$ must be an earlier  move by Abelard in $p_n$. By Lemma \ref{odd strategy}\ (1), since $p_n$ agrees with $S_n$, we have  $\sigma\in(\zeta_0,\zeta_1]$  for some unique pair $\zeta_0\subsetm\zeta_1$ in $S_n$.  If $\sigma$ is not minimal, then $\sigma\subset\zeta_1$ and so Abelard can actually make a move in $p_n$, contrary to the assumption.  Therefore, $\sigma$ is a minimal odd node, and $\vp_\sigma$ is an earlier move by Eloisa in $p_n$. It then follows by  Lemma \ref{odd strategy}\ (2) that the unique even $\tau$ such that $\tau\subsetm\sigma$ belongs to $S_n$. Then there is $m\in V(n)$ such that $\tau=\sigma_m$. The proof that $m$ is unique is as before. 
	
	In both cases, the verification that (a)--(h) hold is straightforward. Note that in case D.2.2  Abelard follows $S_{n+1}$ in $p_{n+1}$ and, by construction, Abelard cannot move in $p_n$, and so (h)  holds.
\end{proof}

\section{Descent recursion} \label{descent recursion}

Fix an elementary recursive pairing function $(\cdot,\cdot)$.  

\begin{notation}
Let $\on$ be an elementary recursive ordinal notation system for some $\varepsilon$-number. We use $\alpha,\beta,\ldots$ to denote elements of $\on$. We rely on context and denote by $<$  the ordering relation on $\on$. By $\on_\alpha$ we denote the set $\{\beta\in\on\mid \beta<\alpha\}$.
\end{notation}
For a precise definition cf.\ \cite[Definition 1.1]{FS95}. Roughly put, $\on$ is what we expect it to be: an elementary recursive well-ordering on $\N$ closed under elementary recursive operations $\alpha+\beta$, $\alpha\cdot \beta$, $\omega^\alpha$, and  satisfying a few basic  properties from ordinal arithmetic.  We will be working with a notion of descent recursion which is  slightly more liberal  than Friedman and Sheard's original definition \cite[Definition 1.5]{FS95}. However, it captures the same class of functions (cf.\ \cite[Proposition 1.9]{FS95}).

\begin{definition}[descent recursion]\label{descent}
Let $\alpha\in\on$,  $g_0\colon \N\to\N\times \on_\alpha$,  $g_1\colon \N^2\times\on\to \N\times\on$ and $g_2\colon \N\to \N$. Set 
\[ g(x)=D_\alpha(g_0,g_1,g_2)(x)=y,\]
where the  sequence $(y_0,\alpha_{0})(y_{1},\alpha_{1})\cdots (y_{s},\alpha_{s})$ is defined by 
\begin{align*}
	g_0(x)& = (y_{0},\alpha_{0}), \\
	g_1(x,y_i,\alpha_i) &= (y_{i+1},\alpha_{i+1}),
\end{align*}
$s$ is least such that $\alpha_s\leq\alpha_{s+1}$, and the output $y=g_2(y_s)$.  We say that $g$ is obtained from the $g_i$'s by a descent recursion below $\alpha$. The definition extends  to functions with more than one argument.
\end{definition}

From now on, by elementary we mean elementary recursive. 
\begin{definition}
The class $\DR\restr\on$ of descent recursive functions on $\on$ is defined by:	
\begin{itemize}[leftmargin=5mm]
	\item elementary  functions are in $\DR\restr \on$;
	\item if $g_0,g_1,g_2$ are in $\DR\restr\on$ and $\alpha\in\on$, then $D_\alpha(g_0,g_1,g_2)$ is in $\DR\restr\on$.
\end{itemize}
\end{definition}
\begin{definition}
A function   $g$ is  $\alpha$-descent recursive ($\alpha$-$\DR$) if $g= D_\alpha(g_0,g_1,g_2)$, where $g_0$, $g_1$, $g_2$   are elementary.
\end{definition}
\begin{proposition}
Every function  in $\DR\restr\on$ is $\alpha$-$\DR$ for some $\alpha$.	
\end{proposition}
\begin{proof}
Elementary functions are $\alpha$-$\DR$ for every $\alpha>0$. Suppose $g=D_\alpha(g_0,g_1,g_2)$ is obtained  from  $\alpha_i$-$\DR$ functions $g_i$ for $i<3$.   To compute  $g(x)$, we first compute $g_0(x)$ in $\alpha_0$ elementary steps. Then we iterate $g_1$ $\alpha$-many times. Each iteration requires $\alpha_1$ elementary steps. Finally, we compute the output using $g_2$ in $\alpha_2$ elementary steps. In total we need $\alpha_2+\alpha_1\alpha+\alpha_0$ elementary steps. 
\end{proof}

\begin{proposition}
The class $\DR\restr\on$ is closed under composition.	
\end{proposition}
\begin{proof}
If $h$ is $\alpha$-$\DR$ and $g$ is $\beta$-$\DR$, then $g\circ h$ is $(\beta+\alpha)$-$\DR$.	
\end{proof}

\begin{example}[primitive recursion]
\begin{enumerate}[leftmargin=0.8cm]
\item The functions definable (from elementary  functions) with one application of primitive recursion are $\omega$-$\DR$.
\item The functions definable (from elementary  functions) with $n$ nested applications of primitive recursion are $\omega^n$-$\DR$.
	\item $\DR\restr\omega^\omega$ is the class of primitive recursive functions.
\end{enumerate}

Note that the totality of an $\omega^n$-$\DR$ function can be proved in ${\sf I}\Sigma_1$ (primitive recursive well-foundedness of $\omega^n$). On the other hand, it is well-known (Parsons) that the provably recursive functions of ${\sf I}\Sigma_1$  are exactly the primitive recursive ones. 
\end{example}

For later use, we give the following definition.

\begin{definition}[stages of computation]
	Let  $g\colon\N\to \N$ be $\alpha$-$\DR$. Let $x\in\N$ and $(y_0,\alpha_0)(y_{1},\alpha_{1})\cdots (y_{s},\alpha_{s})$ be the computation of $g(x)$. The {\em initial} stage of $g(x)$ is $(y_0,\alpha_0)$. The {\em final}  stage of $g(x)$ is $(y_s,\alpha_s)$.
	Every stage $(y_i,\alpha_i)$ before the final stage is an  {\em internal} stage of $g(x)$.
\end{definition}

\section{Games for Peano Arithmetic and complexity} \label{peano games}

We  assume that our language has symbols for sufficiently many elementary  recursive functions and $=$ is the only relation symbol. We then assume an elementary recursive  coding of terms, formulas and finite sequences. We may safely assume that the collection of closed true literals is primitive recursive. However, this will not play any role.  In fact, when considering winning strategies for theorems of $\pa$, we will just use the fact that for any finite set $F$ of formulas and terms, it is elementary to decide whether a closed literal is true (to compute the numerical value of a closed term), as long as such literal (term) is a closed instantiation with numerals of some literal (term) from the set $F$. 

\begin{definition}[translation]
	\begin{align*}
		\vp\lor\psi &\simeq \E \{\vp,\psi\} &  \vp\land\psi &\simeq \A\{\vp,\psi\} \\
		\exists x\, \vp(x) &\simeq \E \{ \vp(e) \mid e \text{ closed term}\} & \forall x\, \vp(x) &\simeq \A_{n\in\N} \vp(\bar n) 
	\end{align*}	
\end{definition}

To make sure that we only use elementary functions, we introduce the following notion.   

\begin{definition}
	A strategy $f$ for Eloisa in $G(\Gamma)$ is	{\em finitely guessing} if there is a finite set $F$ of terms  such that  whenever a move of the form $a= \guess \vp_x(e)$ appears in $p\in T(f)$, then $e$   is a  closed instantiation with numerals of some term from  $F$.
\end{definition}

\begin{remark}
	If $f$ is finitely guessing, then  we have an  elementary  function that on input $e$ outputs its numerical value, whenever the closed term $e$ occurs in some play $p\in T(f)$. In particular, we  have an elementary  procedure to check the truth of closed literals  whenever they appear within a play $p\in T(f)$.\footnote{In general, without any restriction,  it would take $\omega$ elementary steps (i.e., an $\omega$-$\DR$ function) to evaluate a closed term.  Similarly,   if we had (symbols for) all primitive recursive functions,   it would take in general  $\omega^\omega$ elementary steps to compute the value of a closed term.}  
\end{remark}

\begin{convention}
	From now on, let $\on$ be an elementary recursive ordinal notation system for some $\varepsilon$-number. 
\end{convention}

By $g=(f,h)$ we denote the function $g(x)=(f(x),h(x))$.

\begin{definition}[descent recursive strategy]
Let $\gamma,\alpha\in \on$ and  $f,\hh\colon\N^{<\N}\to \N$. We say that $g=(f,\hh)$ is a $(\gamma,\alpha)$-strategy for Eloisa in $G(\Gamma)$ if the following hold:
\begin{itemize}[leftmargin=5mm]
		\item $g$ is $\gamma$-$\DR$;
		\item $f$ is a  non-losing strategy for Eloisa in $G(\Gamma)$;
		\item $\hh$ is a  height function  on $T(f)$ with  $\hh(\lr)=\alpha$.
\end{itemize} 
When $\gamma=\alpha$, we say that $g$ is a $\gamma$-strategy. A descent recursive strategy $(f, \hh)$ is finitely guessing if so is $f$. 
\end{definition}

\begin{remark}
A $(\gamma,\alpha)$-strategy is winning.  
\end{remark}

\begin{theorem}[cut elimination]\label{cut}
Let $\vp$ be an $\E$-formula. Let $g_0=(f_0,\hh_0)$ and $g_1=(f_1,\hh_1)$ be finitely guessing $(\gamma,\alpha)$-strategies for Eloisa  in $G(\Gamma,\neg\vp)$ and in $G(\Gamma,\vp)$, respectively. Then there is  a finitely guessing $c_\nu(\eta)$-strategy  $g=(f,\hh)$ for Eloisa in $G(\Gamma)$, where  $\eta=\gamma\cdot (\alpha+2)$, and   $\nu$ is a bound on the depth of $\vp$.
\end{theorem}

\begin{proof}
Suppose $\vp$ has depth $\leq \nu$.  The proof of Theorem \ref{height} gives us  a primitive recursive function $c(\cdot,\cdot,\cdot)$ such that  $c(\nu,\alpha,\cdot)$ is a height function on $T=T(\nu,\alpha)$ so that  $\height(T)\leq c_\nu(\alpha)$.   For any fixed $\nu$, the function $c(\nu,\cdot,\cdot)$ is elementary. From now on, let $c(\cdot)=c(\nu,\eta,\cdot)$. We compute the cut strategy $f=\cut(f_0,f_1)$  by using the canonical debate (Definition \ref{canonical debate}). Note that every step in the canonical debate for a sequence $p$ is elementary recursive (relative to $f_0$ and $f_1$). 
	
Observe that the cut strategy inherits the property of being finitely guessing from $f_0$ and $f_1$ by taking $F=F_0\cup F_1$, where $F_i$ witnesses the property for $f_i$.  

We define a $c_\nu(\eta)$-$\DR$ function by showing how to compute in an elementary way the initial stage and the next stage of a computation  while descending below $c_\nu(\eta)$.  We define  $\hh(p)=\delta$ whenever $(z,\delta)$ is the final stage of the computation. In particular, we just need to specify $f(p)$ when we reach such final stage.

Suppose $p=a_1\cdots a_r$ is a play in $G(\Gamma)$. The initial stage of $g(p)$ is $(z,\delta)$, where:
 
\begin{itemize}[leftmargin=5mm]
\item 	$z$ consists of:
\begin{itemize}[leftmargin=2mm] 
\item the initial round $d_0$ of the  debate for $p$ (recall that $d_0=p_0$ with $p_0=\lr$);
\item the initial stage $(y_0,\gamma_0)$ of $g_0(p_0)$;
\item the ordinal interaction sequence $u=\beta_0$ with $\beta_0=\gamma(\alpha+1)+\gamma_0$; 
\end{itemize}
\item $\delta=c(u)$.  
\end{itemize}
Note that $c(u)< c(\lr)=c_\nu(\eta)$.   In general, any  stage of $g(p)$  is of the form $(z,\delta)$, where 
\begin{itemize}[leftmargin=5mm]
\item   $z$  consists of:
\begin{itemize}[leftmargin=2mm] 
\item  some round  $d_s=e_{s}\, p_n= (p_0,\sigma_0, i_0)\cdots (p_{n-1},\sigma_{n-1}, i_{n-1})\,  p_{n}$ of a debate for $p$ with $\pGamma_{s}=a_1\cdots a_l\subseteq  p$;
\item some stage  $(y_m,\gamma_m)$ of $g_{\party(m)}(p_m)$  for each $m\leq n$ which is final for $m<n$;
\item an ordinal interaction sequence $u=  \hat{u}\,  \beta_n=(\beta_0,i_0)(\beta_1,i_1)\cdots(\beta_{n-1},i_{n-1})\, \beta_n$ in $T(\nu,\eta)$, where $\beta_m$ is of the form $$\gamma\alpha_m+\gamma_m  \text{ for every $m<n$ with } \alpha_m=\hh_{\party(m)}(p_m) \text{ and } \gamma_m<\gamma,$$  and $\beta_n$ is of the form $\gamma\rho_n+\gamma_n$ with $\gamma_n<\gamma$;
\end{itemize}
 \item $\delta=c(u)$.
\end{itemize}
Recall that $\party(n)$ is the parity of $n$. \smallskip

We run through all cases of the canonical debate (Definition \ref{canonical debate}) by taking the internal stages of the two strategies $g_0$ and $g_1$ into further account. To ease notation, we rely on context and write $g_i$, $f_i$, $\hh_i$ instead of $g_{\party(m)}$, $f_{\party(m)}$, $\hh_{\party(m)}$.  \smallskip

Construction:\smallskip

If $(y_n,\gamma_n)$ is an internal stage of $g_i(p_n)$,  then $(z,\delta)$ is an internal stage of $g(p)$. Update $(y_n,\gamma_n)$ to $(y'_n,\gamma'_n)=g_{i,1}(p_n,y_n,\gamma_n)$ and  $u$ to $u'=\hat{u}\, (\beta_n')$, where $\beta_n'=\gamma\rho_n+\gamma'_n$. Note that $\delta>c(u')=\delta'$. (Here,  $g_{i,1}$ is the   elementary function   that computes the next stage  in the computation of $g_i$.)\smallskip

Otherwise, $(y_n,\gamma_n)$ is the final stage of $g_i(p_n)$. Say $f_i(p_n)=a$ and $\hh_i(p_n)=\alpha_n$.  \smallskip

Case 1 (cf.\  A). The play $a_1\cdots a_l$ is winning. Then  $(z,\delta)$ is the final stage of $g(p)$ and  $f(p)=a$. \smallskip

Case 2 (cf.\  B.1 and D.2.2).  Abelard  responds to a query on the $\vp$ side in $p_{i_n}$.  Say Abelard plays  $b$ so that $p_{n+1}=p_{i_n}b$. Then $(z,\delta)$ is an internal stage of $g(p)$. Define $(y_{n+1},\gamma_{n+1})$ to be the first stage of $g_i(p_{n+1})$. Update $d_{s}$ to   $e_s\, (p_n,\sigma_n,i_n)\, p_{n+1}$ and $u$ to  $u'=\hat{u}\, (\beta'_n,i_n)\, \beta_{n+1}$, where
\begin{itemize}[leftmargin=5mm]
\item  $\beta'_n=\gamma\alpha_n+\gamma_n$, 
\item  $\beta_{n+1}=\gamma\alpha_{i_n}+\gamma_{n+1}$.
\end{itemize} \smallskip

Case 3. It is Eloisa's turn in $p_n$. 

Case 3.1 (cf.\  C.1.1).  $a$ is a move on the $\Gamma$ side and  either  $l=r$ or $l<r$ with $a\neq a_{l+1}$. Then  $(z,\delta)$ is the final stage of $g(p)$ and  $f(p)=a$.

Case 3.2 (cf.\  C.1.2).  $a$ is a move on the $\Gamma$ side,   $l<r$ and $a=a_{l+1}$.  Then $(z,\delta)$ is an internal stage. Update $d_{s}$ to $d_{s+1}=e_{s}\, (p_na)$. Update  $(y_n,\gamma_n)$ to  the initial stage $(y_n',\gamma_n')$ of $g_i(p_na)$.
Update $u$ to $u'=\hat{u}\, \beta_n'$, where $\beta_n'=\gamma\alpha_n+\gamma_n'$.\smallskip 

Case 3.3 (cf.\  C.2) $a$ is a move on the $\vp$ side. Proceed as in case 3.2. \smallskip

Case 4.   It is Abelard's turn in $p_n$ .

Case  4.1 (cf.\  D.1.1). The last move in $p_n$ is a query on the $\Gamma$ side and $l=r$ or $l<r$ and $a_{l+1}$ is not a legal move.  Then  $(z,\delta)$ is the final stage of $g(p)$ and  $f(p)=a$. 

Case  4.2 (cf.\ D.1.2). The last move in $p_n$ is a query on the $\Gamma$ side  with $l<r$ and $b=a_{l+1}$ is a legal move. Then $(z,\delta)$ is an internal stage. Update $d_{s}$ to $d_{s+1}=e_{s}\, (p_nb)$,  the stage $(y_n,\gamma_n)$ to the initial stage $(y_n',\gamma_n')$ of $g_i(p_nb)$, and $u$ to $u'=\hat{u}\, \beta_n'$, where $\beta_n'=\gamma\alpha_n+\gamma_n'$. 

Case 4.3 (cf.\ D.2.1.)  The last move in $p_n$ is a query on the $\vp$ side and Abelard  can move in $p_n$. Say Abelard plays $b$.  Proceed as in case 4.2.\smallskip

End of construction. \smallskip

One can verify by induction that at any stage of the computation $\alpha_m=\hh_i(p_m)$ for every $m<n$ and  $\rho_n>\alpha_n=\hh_i(p_n)$. It is then easy to check that $f=\cut(f_0,f_1)$ and $t$ is a height function on $T(f)$. The details are left to the reader. 
\end{proof}

\begin{theorem}\label{sound}
	Every theorem of $\pa$ has a  finitely guessing $(\gamma,\alpha)$-strategy for some $\gamma,\alpha<\varepsilon_0$.
\end{theorem}

\begin{proof}
By induction on a proof $\pi(x_1,\ldots,x_k)$ of  $\Gamma$ in $\pa$, we show that for some $\gamma,\alpha<\varepsilon_0$ there is a $\gamma$-$\DR$  function $g=(f,h)$ with parameters, namely $g\colon \N^k\times \N^{<\N}\to \N$,  such that $g(n_1,\ldots,n_k,\cdot )$ is  a $(\gamma,\alpha)$-strategy for Eloisa in the game  $\Gamma(\bar n_1,\ldots,\bar n_k)$ for every $n_1,\ldots, n_k\in\N$. Moreover, for every play $pa\in T(f)$, where $a$ is a move of the form $\guess \vp_x(e)$,  $e$   is a  closed instantiation with numerals of some term appearing in $\pi$. In particular, $f(n_1,\ldots,n_k,\cdot )$ is finitely guessing.   
	
	\begin{claim}
		Basic axioms 
	\end{claim}
	If $\Gamma$ is a basic axiom, then  every elementary function $(f,\hh)$ with $\hh(\lr)=0$ can be regarded as a $(1,0)$-strategy for Eloisa.
	
	\begin{claim}
		Normal rules
	\end{claim}
	Let us treat $\exists$ and $\forall$ only. 
	
	($\exists$-rule) Suppose $\pi(x_1,\ldots,x_k)\vdash \Gamma, \exists x\, \vp(x)$ with $\pi_0(x_1,\ldots,x_k)\vdash \Gamma, \vp_x(e)$.  By induction, let $g_0(n_1,\ldots,n_k,\cdot)$ be a $(\gamma,\alpha)$-strategy in the game  $$\Gamma(\bar n_1,\ldots,\bar n_k),\vp_x(e)_{x_1,\ldots,x_k}(\bar n_1,\ldots,\bar n_k).$$ Let $f$ be defined by 
	\begin{align*}
		f(n_1,\ldots,n_k,\lr) &= \vp_x(e)_{x_1,\ldots,x_k}(\bar n_1,\ldots,\bar n_k); \\
		f(n_1,\ldots,n_k,a p) & =f_0(n_1,\ldots,n_k,p).
	\end{align*}
    Let 
    \begin{align*}
    	\hh(n_1,\ldots,n_k,\lr) &= \alpha+1; \\
    	\hh(n_1,\ldots,n_k,a p) & =\hh_0(n_1,\ldots,n_k,p).
    \end{align*}
	Then $(f,\hh)(n_1,\ldots,n_k,\cdot)$ is a $(\gamma,\alpha+1)$-strategy in the game $\Gamma(\bar n_1,\ldots,\bar n_k), \exists x\, \vp(\bar n_1,\ldots,\bar n_k,x)$.
	\smallskip

	($\forall$-rule) Suppose $\pi(x_1,\ldots,x_k)\vdash \Gamma, \forall  x\, \vp(x)$ with $\pi_0(x_1,\ldots,x_k,y)\vdash \Gamma, \vp_x(y)$. By induction, suppose $g_0(n_1,\ldots,n_k,n,\cdot)$ is a $(\gamma,\alpha)$-strategy in the game $$\Gamma(\bar n_1,\ldots,\bar n_k), \vp(\bar n_1,\ldots,\bar n_k,\bar n).$$  Let $f(n_1,\ldots,n_k,\lr)=a$, where $a$ is  $\query \forall x\, \vp(\bar n_1,\ldots,\bar n_k,x)$, and $f(n_1,\ldots,n_k,a b p)=f_0(n_1,\ldots,n_k,n,p)$ whenever  $b$ is   $\reply \vp(\bar n_1,\ldots,\bar n_k,\bar n)$. 
	Then $(f,\hh)(n_1,\ldots,n_k,\cdot)$ is a $(\gamma,\alpha+2)$-strategy in the game $\Gamma(\bar n_1,\ldots,\bar n_k), \forall  x\, \vp(\bar n_1,\ldots,\bar n_k,x)$ by choosing the obvious  $h$.
	
\begin{claim}
		Induction
\end{claim}
Eloisa has a $(1,\omega+1)$-strategy in the game \[   \neg\vp(0), \exists x\, (\vp(x)\land  \neg\vp(x+1)), \forall x\, \vp(x). \] 
Let $c<\omega$ be the rank of $\vp$.   We  display the strategy tree in Figure \ref{figure uno} (for ease of notation, we identify the numeral $\bar n$ with the number $n$).    We make no use of internal stages save for the initial one (no cuts are involved). It is clear how to fill the dots in Figure \ref{figure uno}. Note that Eloisa has a winning copycat strategy as soon as Abelard plays $\vp(0)$ or both $\vp(n)$ and $\neg\vp(n)$ with $n>0$,  ending the game in $\leq 3c$ steps. In Figure \ref{figure uno},   $\delta(n)$ is $\vp(n)$ if $\vp(n)$ is an $\A$-formula or  $\neg\vp(n)$ otherwise. If $\vp$ is a literal, then Eloisa wins without further moves.
	
	% forest 1
	\begin{figure} %[h!]
		\begin{forest} for tree={s sep=20mm, inner sep=5}
			[ {$\omega+1$} 
			[ {\color{gray} $\query$}  $\forall x\, \vp(x)$ {\color{gray} at} $\omega$, calign=center
			[ {\color{gray} $\reply $} $\vp(0)$ {\color{gray} at $\cdots$} [{\color{gray} $\query$} $\delta(0)$ {\color{gray} at} $3c$  [$\vdots$]]] 
			[ {\color{gray} $\reply $} $\vp(n+1)$ {\color{gray} at $\cdots$}
			[ {\color{gray} $\guess$}  $\vp(n)\land\neg\vp(n+1)$ {\color{gray} at $\cdots$}
			[ {\color{gray} $\query$}  $\vp(n)\land\neg\vp(n+1)$ {\color{gray} at $\cdots$}, 
			calign=last
			[{\color{gray} $\reply $} $\neg\vp(n+1)$ {\color{gray} at $\cdots$} [{\color{gray} $\query$} $\delta(n+1)$ {\color{gray} at} $3c$ [$\vdots$ ] ]] [ {\color{gray} $\reply $} $\vp(n)$ {\color{gray} at} $3c+3n+2$
			 [ {\color{gray} $\guess$} $\vp(n-1)\land\neg\vp(n)$ {\color{gray} at} $3c+3n+1$ [ {\color{gray} $\query$} $\vp(n-1)\land\neg\vp(n)$ {\color{gray} at} $3c+3n$, calign=last
			[ {\color{gray} $\reply $} $\neg\vp(n)$ {\color{gray} at $\cdots$} [{\color{gray} $\query$} $\delta(n)$ {\color{gray} at} $3c$ [$\vdots$]]] [ {\color{gray} $\reply $} $\vp(n-1)$ {\color{gray} at $\cdots$} [$\vdots$] ] ]  ]  ]   ] ]  ] ] ] 
		\end{forest}
		%\vspace{20pt}
		\caption{Strategy tree for  induction.} \label{figure uno}
	\end{figure}
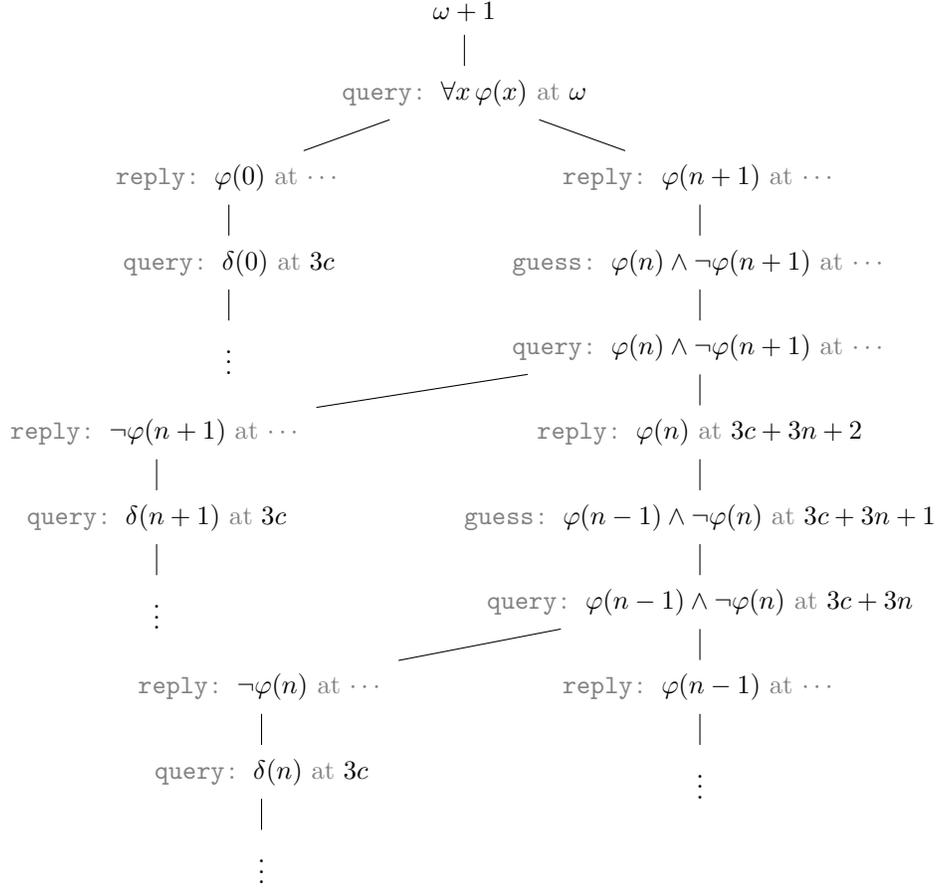 
	
	\begin{claim}
		Cut
	\end{claim}	
	The proof of Theorem \ref{cut} relativizes by allowing parameters.  If the cut formula is a literal, say $\vp(x_1,\ldots,x_k)$, then the  strategy
	\[  g(n_1,\ldots,n_k,p)=\begin{cases}
		g_0(n_1,\ldots,n_k,p) & \text{ if } \vp(\bar n_1,\ldots,\bar n_k) \text{ is true;} \\
		g_1 (n_1,\ldots,n_k,p)  & \text{ otherwise,}
	\end{cases} \]
	is a $(\gamma,\alpha)$-strategy \sloppy in $G(\Gamma(\bar n_1,\ldots,\bar n_k))$, if $g_0(n_1,\ldots,n_k,\cdot)$ is a $(\gamma,\alpha)$-strategy in $G(\Gamma(\bar n_1,\ldots,\bar n_k),\neg\vp(\bar n_1,\ldots,\bar n_k))$ and $g_1(n_1,\ldots,n_k,\cdot)$ is a $(\gamma,\alpha)$-strategy in 
	$G(\Gamma(\bar n_1,\ldots,\bar n_k),\vp(\bar n_1,\ldots,\bar n_k))$. 
	
\end{proof}

\subsection{Transfinite induction}

\begin{definition}[transfinite induction]
	Let $\alpha\in\on$. The schema $\TI(\alpha)$ of transfinite induction up to $\alpha$ consists of 
	\[ \forall \beta\, (\forall\gamma<\beta\, \vp(\gamma)\imp \vp(\beta)) \imp \forall \gamma<\alpha\, \vp(\gamma), \]
	where $\vp$ ranges over all formulas. Let $\TI\restr\on=\bigcup_{\alpha\in\on}\TI(\alpha)$.
\end{definition}

\begin{theorem}\label{transfinite}
	For every $\alpha\in\on$  and for every formula $\vp$ there is a finitely guessing $(1,3c+5\alpha+2)$-strategy for Eloisa in the game 
	\[   \exists \beta\, (\forall \gamma<\beta\, \vp(\gamma)\land \neg\vp(\beta)), \forall \gamma<\alpha\, \vp(\gamma) \] for some $c<\omega$.
\end{theorem}
\begin{proof}
	Let $c$ be the rank of $\vp$. We display the strategy tree in Figure \ref{figure due}. Note that as soon as Abelard plays $\neg\vp(\beta)$, Eloisa has a  winning copycat strategy that ends the game in $\leq 3c$ steps. For example, suppose $\vp(\beta)=\forall x\, \exists y\, \vartheta_\beta(x,y)$, where $\vartheta_\beta(x,y)$ is a literal. Then Eloisa wins  by making the  moves in Figure \ref{figure tre} (again, we write $n$ in place of the numeral $\bar n$).
	
	% forest 2
	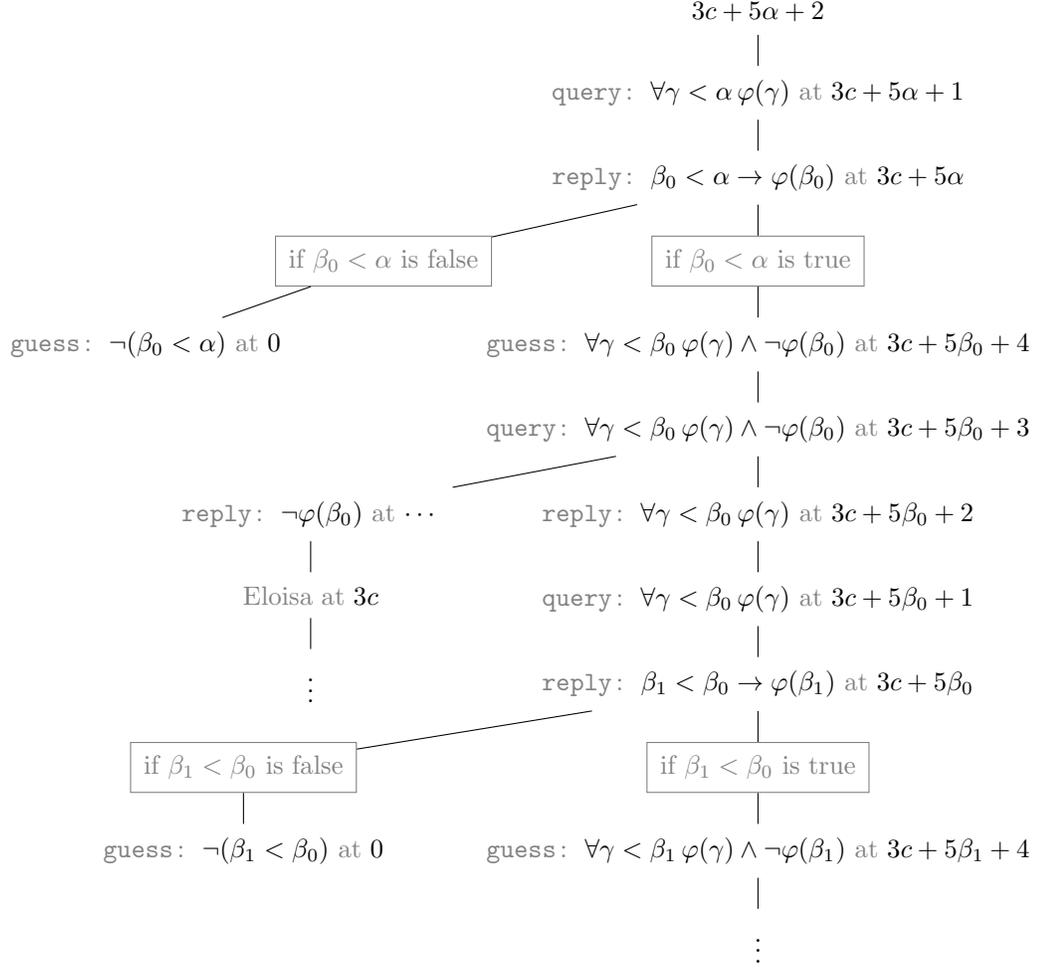
\begin{figure} %[h!]
		\begin{forest}  for tree={s sep=10mm, inner sep=5, l=-5}
			[ {$3c+5\alpha+2$} % root
			[ {\color{gray} $\query$} $\forall \gamma<\alpha\, \vp(\gamma)$ {\color{gray} at}  $3c+5\alpha+1$ % child
			[ {\color{gray} $\reply $} $\beta_0<\alpha\imp \vp(\beta_0)$ {\color{gray} at} $3c+5\alpha$,  calign=last   % child 
			[ if $\beta_0<\alpha$ is false, draw, gray, calign=last
			[ {\color{gray} $\guess$} $\neg(\beta_0<\alpha)$ {\color{gray} at} $0$ ] [, phantom] ]  % left
			[ if $\beta_0<\alpha$ is true, draw, gray % right
			[ {\color{gray} $\guess$}  $\forall \gamma<\beta_0\, \vp(\gamma)\land \neg\vp(\beta_0)$ {\color{gray} at} $3c+5\beta_0+4$  % child
			[ {\color{gray} $\query$}  $\forall \gamma<\beta_0\, \vp(\gamma)\land \neg\vp(\beta_0)$ {\color{gray} at} $3c+5\beta_0+3$, calign=last % child
			[ {\color{gray} $\reply $}  $\neg\vp(\beta_0)$ {\color{gray} at} $\cdots$  [{\color{gray} Eloisa at} $3c$ % left
			[$\vdots$]  ]  ] 
			[  {\color{gray} $\reply $} $\forall \gamma<\beta_0\, \vp(\gamma)$ {\color{gray} at} $3c+5\beta_0+2$ % right
			[ {\color{gray} $\query$} $\forall \gamma<\beta_0\, \vp(\gamma)$  {\color{gray} at}  $3c+5\beta_0+1$ % child
			[ {\color{gray} $\reply $} $\beta_1<\beta_0\imp \vp(\beta_1)$ {\color{gray} at}  $3c+5\beta_0$,  calign=last   % child
			[ if $\beta_1<\beta_0$ is false, draw, gray  %left
			[ {\color{gray} $\guess$} $\neg(\beta_1<\beta_0)$ {\color{gray} at} $0$  ]  ] % child
			[ if $\beta_1<\beta_0$ is true, draw, gray   % right
			[ {\color{gray} $\guess$}  $\forall \gamma<\beta_1\, \vp(\gamma)\land \neg\vp(\beta_1)$ {\color{gray} at} $3c+5\beta_1+4$    % child
			[$\vdots$ ] ]  
			]  ]  ]   ] ]    
			] ]  ]   ]  ]
		\end{forest}
		\caption{Strategy tree for transfinite induction up to $\alpha$.} \label{figure due}
	\end{figure}
	
	% forest 3 
	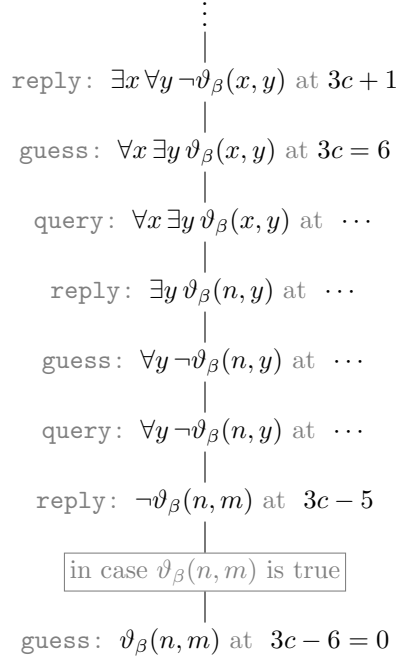
\begin{figure} [!h]
		\begin{forest} for tree={s sep=5mm, inner sep=2, l=0}
			[ $\vdots$
			[ {{\color{gray} $\reply $} $\exists x\, \forall y\, \neg\vartheta_\beta(x,y)$ {\color{gray} at} $3c+1$} 
			[ {{\color{gray} $\guess$} $\forall  x\, \exists y\,  \vartheta_\beta(x,y)$ {\color{gray} at} $3c=6$}  	
			[ {{\color{gray} $\query$} $\forall  x\, \exists y\,  \vartheta_\beta(x,y)$ {\color{gray} at } $\cdots$ } 
			[ {{\color{gray} $\reply $} $\exists y\, \vartheta_\beta( n,y)$ {\color{gray} at } $\cdots$}   
			[ {{\color{gray} $\guess$} $\forall y\, \neg \vartheta_\beta( n,y)$ {\color{gray} at } $\cdots$ }  
			[ {{\color{gray} $\query$} $\forall y\, \neg \vartheta_\beta( n,y)$ {\color{gray} at } $\cdots$ } 
			[ {{\color{gray} $\reply $} $\neg \vartheta_\beta( n, m)$ {\color{gray} at } $3c-5$ } 
			[ {in case $\vartheta_\beta(n, m)$ is true}, draw, gray, 
			[ {{\color{gray} $\guess$} $\vartheta_\beta( n, m)$ {\color{gray} at } $3c-6=0$ }          
			]   ]   ]   ]   ]  ] ] ] ] ]
		\end{forest}
		\caption{Example of copycat strategy for $\forall\exists$ formulas.} \label{figure tre}
	\end{figure}
\end{proof}

\begin{theorem}
	Every theorem of $\pa+\TI\restr\on$ has a finitely guessing $(\gamma,\alpha)$-strategy for some $\gamma,\alpha\in\on$.
\end{theorem}
\begin{proof}
	The proof of Theorem \ref{sound} generalizes in virtue of  Theorems \ref{cut} and  \ref{transfinite}.
\end{proof}

\subsection{Provably recursive functions}

\begin{theorem}\label{provably}
	Let $\vp=\forall x\, \exists y\, \vp_0(x,y)$ be a $\Pi_2$ sentence, where $\vp_0(x,y)=\exists z\, \vartheta(x,y,z )$  with $\vartheta(x,y,z)$ elementary. If there is a finitely guessing $(\gamma,\alpha)$-strategy $g=(f,\hh)$ for Eloisa in $G(\vp)$, then there is a  $\gamma\cdot(\alpha+2)$-$\DR$ function $\hat g\colon \N\to \N$ such that   $\vp_0(n,\hat g(n))$ for every $n$. 
\end{theorem}
\begin{proof}
	To compute $\hat g(n)$ we play $G(\vp)$ according to $f$ and let Abelard  respond to  $q=\query \vp$ with $b=\reply \exists y\, \vp_0(\bar n,y)$.  Eloisa is  forced   to either play $q$ or to  make a guess of the form  $\vp_0(\bar n,e_0)$ or $\vartheta(\bar n,e_0,e_1)$ for some closed terms $e_0$ and $e_1$. Eventually Eloisa will make a winning guess $\vartheta(\bar n, e_0,e_1)$ so that $\vp_0(n,m)$ holds where $m$ is the value of $e_0$. Note that we can evaluate any such closed term in an elementary way since $f$ is finitely guessing. Essentially, we are looking for the unique winning play $p=a_1\cdots a_m\in T(f)$ such that $a_{i+1}=b$ whenever $a_i=q$ for every $i<m$.  
	That all this can be done via a  descent recursion below $\gamma(\alpha+2)$ is clear.  We start the computation of $\hat g(n)$ at   $(y_0,\gamma(\alpha+1)+\gamma_0)$, where $(y_0,\gamma_0)$ is the initial stage of $g(\lr)$.  
	We then compute the next stage  or stop the computation according  to the following procedure. If $(y_s,\gamma_s)$ is not the final stage of $g(p)$, we are in an internal stage of $\hat g(n)$. The next stage  is $(y_{s+1},\gamma\alpha_s+\gamma_{s+1})$, where $(y_{s+1},\gamma_{s+1})$ is the next stage of $g(p)$. 
	If $(y_s,\gamma_s)$ is the final stage of $g(p)$, say $f(p)=a$ and $\hh(p)=\alpha_0$, we act according to whether $a$ is the query $q$ or some guess.   In the first case, the next stage of $\hat g(n)$ is $(y_0,\gamma\alpha_0+\gamma_0)$, where $(y_0,\gamma_0)$ is the initial stage of $g(pqb)$. In the second case, if $a$ is not a winning guess,  the next stage of $\hat g(n)$ is $(y_0,\gamma\alpha_0+\gamma_0)$, where  $(y_0,\gamma_0)$ is the initial stage of $g(pa)$. Finally, if $a$ is a winning guess, say $\vartheta(\bar n, e_0,e_1)$,  we evaluate  $e_0$ and output $\hat g(n)=\text{ the value of } e_0$. 
\end{proof}

\begin{corollary}[provably recursive functions]
	$\DR\restr\varepsilon_0$ is the class of provably recursive functions of $\pa$. More in general, $\DR\restr\on$ is the class of provably recursive functions of $\pa+\TI\restr\on$. 
\end{corollary}
\begin{proof}
	One direction follows from (the proof of) Theorem \ref{sound} and  Theorem \ref{provably}. Conversely, we use $\TI(\alpha)$ to prove that every given $\alpha$-$\DR$ function is total. 
\end{proof}

\subsection{No-counter-example}

We start by lifting Definition \ref{descent} to type $2$ functionals  in the obvious way.

\begin{definition}[descent recursive functionals]
	A functional   $G\colon \N^\N\to \N$ is $\alpha$-$\DR$ if $G=D_\alpha(G_0,G_1,G_2)$, where $G_0\colon \N^\N\to \N\times\on_\alpha$, $G_1\colon \N^\N\times\N\times \on\to \N\times\on$ and $G_2\colon\N\to\N$ are elementary recursive. 	The definition extends to functionals with more than one argument in the usual way.  
\end{definition}

\begin{definition}
	Let $\vp$ be  $\exists x_1\, \forall y_1\, \cdots\,  \exists x_k\, \forall y_k\, \vartheta(x_1,y_1,\ldots,x_k,y_k)$.  Let $G=(G_1,\ldots,G_k)$ be a tuple of functionals 
	\[ G_j\colon \overbrace{\N^{\N}\times \cdots \times \N^\N}^{k\text{ times}}\to \N, \ \quad j\in[1,k]. \]  We say that $G$ satisfies the {\em no-counter-example interpretation} of $\vp$, in short $G \text{ n.c.i. } \vp$, if 
	\[ \forall f_1\cdots \forall f_k\, \vartheta(x_1,f_1(x_1),\ldots, x_k, f_k(x_1,\ldots,x_k)),  \]
	with $x_j=G_j(f_1,\ldots,f_k)$ for $j\in[1,k]$. 
\end{definition}

\begin{theorem}
	Let $\vp$ be $\exists x_1\, \forall y_1\, \cdots\,  \exists x_k\, \forall y_k\, \vartheta(x_1,y_1,\ldots,x_k,y_k)$ with $\vartheta(x_1,\ldots,y_k)$ elementary. If there is a finitely guessing $(\gamma,\alpha)$-strategy $g=(f,\hh)$ for Eloisa in $G(\vp)$, then there are  $\gamma\cdot(\alpha+2)$-$\DR$ functionals $G$ such that $G\text{ n.c.i. }\vp$.  
\end{theorem}
\begin{proof}[Proof idea]
	Suppose $g=D_\gamma(g_0,g_1,g_2)$. To compute $G(f_1,\ldots,f_k)$ we play $G(\vp)$ according to the strategy $f$. The functions $f_1,\ldots,f_k$  provide a strategy for Abelard. To  a query 
	\[    \forall y_j\, \exists x_{j+1}\, \forall y_{j+1}\, \cdots\,  \exists x_k\, \forall y_k\,   \vartheta(e_1,\bar m_1, \ldots, e_j,y_j,\ldots, x_k,y_k), \]
	Abelard  responds with
	\[    \exists x_{j+1}\, \forall y_{j+1}\, \cdots\, \exists x_k\, \forall y_k\,     \vartheta(e_1,\bar m_1, \ldots, e_j,\bar m_j,\ldots, x_k,y_k), \]
	where $m_j$ is the value of $f_j(n_1,\ldots, n_j)$ and each $n_i$ is the value of $e_i$. Eventually, Abelard  will  make a winning move
	\[  \vartheta(e_1,\bar m_1,,\ldots, e_k,\bar m_k), \]
	so that  $G_j(f_1,\ldots,f_k)=n_j$, where $n_j$ is the value of $e_j$.  Essentially, we are looking for the unique winning play $p\in T(f)$ that agrees with both strategies.  This can be done through a descent recursion below  $\gamma(\alpha+2)$. 
	We start computing $G_j(f_1,\ldots,f_k)$ at $(y_0,\gamma(\alpha+1)+\gamma_0)$, where $(y_0,\gamma_0)$ is the initial stage of $g(\lr)$.  The transition from one stage to the next is defined in a similar manner as in  the proof of Theorem \ref{provably}. Whenever $f(p)$ is a query $a$ on a formula with closed terms $e_1,\ldots,e_j$ and numerical values $n_1,\ldots,n_j$, we obtain the reply $b$ by  computing  $f_j(n_1,\ldots,n_j)$,  and consider the play $pab$ at the next stage.   This is the only place where the $f_j$'s enter the computation. The details are left to the reader. 
\end{proof}

\begin{corollary}[no-counter-example interpretation functionals] 
	Let $\vp$ be a prenex sentence. If $\vp$ is provable in $\pa$ (resp.\ $\pa+\TI\restr\on$), then there is a tuple $G$ of $\alpha$-$\DR$ functionals such that $G \text{ n.c.i. } \vp$, for some $\alpha<\varepsilon_0$ (resp.\ $\alpha\in \on$).  
\end{corollary}

\end{document}